

 \magnification1000 \def\oui{oui}
\ifx\textures\oui
 \input iten.tex
\input CrayolaColors 
 \long\def\rge#1{\Red#1\Black}

 \else


  %
  %

  %
  %

  \catcode`@=12 

 \def\defrefnote#1{\definexref{#1}{{\the\footnotenumber}}{refnotes}}

  %
  %


 \input eplain.tex
 
\ifx\couleurs\oui
\beginpackages
\usepackage{color}
\endpackages 
\long\def\rge#1{{\color{red}#1}}

\fi

\makeatletter
\def\numberedfootnote{%
ÊÊ\global\advance\footnotenumber by 1
ÊÊ\@eplainfootnote{{\number\footnotenumber}}%
}%
\def\makecolumns#1/#2 {\par \begingroup
ÊÊ \@columndepth = #1
ÊÊ \advance\@columndepth by -1
ÊÊ \divide \@columndepth by #2
ÊÊ \advance\@columndepth by 1
ÊÊ \@linestogoincolumn = \@columndepth
ÊÊ \@linestogo = #1
ÊÊ \currentcolumn = 1
ÊÊ \def\@endcolumnactions{%
ÊÊÊÊÊÊ\ifnum \@linestogo<2
ÊÊÊÊÊÊÊÊ \the\crtok \egroup \endgroup \par 
ÊÊÊÊÊÊ\else
ÊÊÊÊÊÊÊÊ \global\advance\@linestogo by -1
ÊÊÊÊÊÊÊÊ \ifnum\@linestogoincolumn<2
ÊÊÊÊÊÊÊÊÊÊÊÊ\global\advance\currentcolumn by 1
ÊÊÊÊÊÊÊÊÊÊÊÊ\global\@linestogoincolumn = \@columndepth
ÊÊÊÊÊÊÊÊÊÊÊÊ\the\crtok
ÊÊÊÊÊÊÊÊ \else
ÊÊÊÊÊÊÊÊÊÊÊÊ&\global\advance\@linestogoincolumn by -1
ÊÊÊÊÊÊÊÊ \fi
ÊÊÊÊÊÊ\fi
ÊÊ }%
ÊÊ \makeactive\^^M
ÊÊ \letreturn \@endcolumnactions
ÊÊ \@columnwidth = \hsize
ÊÊÊÊ \advance\@columnwidth by -\parindent
ÊÊÊÊ \divide\@columnwidth by #2
ÊÊ \penalty\abovecolumnspenalty
ÊÊ \noindent 
ÊÊ \valign\bgroup
ÊÊÊÊ &\hbox to \@columnwidth{\strut \hsize = \@columnwidth ##\hfil}\cr
}%
\makeatother

\lefteqnumbers
   \def\testd{oui}
   \def\choixlat{\ifx\numadroite\testd\righteqnumbers
            \else  \lefteqnumbers\fi}
    \choixlat

\catcode`@=\letter
\def\@eplainfootnote#1{\let\@sf\empty 
  \ifhmode\edef\@sf{\spacefactor\the\spacefactor}\/\fi
  \global\advance\hlfootlabelnumber by 1
  \hlstart@impl{foot}{\hlfootlabel}%
  \hldest@impl{footback}{\hlfootbacklabel}%
  \hbox{$^{(#1)}$}%
  \hlend@impl{foot}%
  \@sf\vfootnote{#1.}%
}%
\catcode`@=\other

  \interfootnoteskip=0pt
  \let\note=\numberedfootnote
  \everyfootnote={\eightpoint\leftskip=5truemm\rightskip5truemm}
  
  \hsize150truemm\vsize 240truemm\hoffset=5truemm
  \def\dimstand{\hsize 150truemm\vsize 240truemm\hoffset=5truemm\voffset=0truemm}
  
  \pretolerance=500\tolerance=1000\brokenpenalty=5000
  \parindent3mm
  
  \countdef\temps=170
  \temps=\time
  \countdef\nminutes=171{\nminutes = \time}
  \countdef\nheure=172
  \def\heure{\begingroup                     
     \temps = \time \divide\temps by 60
     \nheure = \temps                        
     \nminutes = \time
     \multiply\temps by 60
     \advance\nminutes by -\temps            
     \ifnum\nminutes<10 \toks1 = {0}%
     \else\toks1 = {}%
     \fi
     \number\nheure h\the\toks1 \number\nminutes  
  \endgroup}%

  \newcount\chstart
  \chstart=\pageno
 \headline={\ifnum\pageno=\chstart {\hfill} \else {\hss \tenrm --\ \folio\ --\hss}\fi}
  \footline={\hfill}
  \normalbaselines
  \frenchspacing
    \def\dater{\vglue-10mm\rightline{(\the\day/\the\month/\the\year)}}
  \def\dateheure{\vglue-10mm\rightline{(\the\day/\the\month/\the\year,\ \heure)}}

  \newif\ifpagetitre \pagetitretrue
\newtoks\hautpagetitre \hautpagetitre={\hfill}
\newtoks\baspagetitre \baspagetitre={\hfill}
\newtoks\auteurcourant \auteurcourant={\hfill}
\newtoks\titrecourant \titrecourant={\hfill}
\newtoks\hautpagegauche
\newtoks\hautpagedroite
\newtoks\hautpagemilieu
\hautpagemilieu={\tenrm\hfil -- \folio\ -- \hfil}
\hautpagegauche={\ifx\midfolio\oui\the\hautpagemilieu\else\tenrm\folio\hfill\the\auteurcourant\hfill\fi}
\hautpagedroite={\ifx\midfolio\oui\the\hautpagemilieu\else\hfill\the\titrecourant\hfill\tenrm\folio\fi}
\newtoks\baspagegauche \baspagegauche={\hfil}
\newtoks\baspagedroite \baspagedroite={\hfil}
\headline={\ifpagetitre\the\hautpagetitre
\else\ifodd\pageno\the\hautpagedroite\else\the\hautpagegauche\fi\fi }
\footline={\ifpagetitre\the\baspagetitre
\else\ifodd\pageno\the\baspagedroite
\else\the\baspagegauche\fi\fi \global\pagetitrefalse}

\def\pageblanche{\vfill\eject\pagetitretrue
\null\vfill\eject
\pagetitretrue
}
\def\chgtpage{\ifodd\pageno \else
\pageblanche \fi \pagetitretrue\titreun=0\footnotenumber=0}

\def\chgtpageincrtitreun{\ifodd\pageno \else
\pageblanche \fi \pagetitretrue\footnotenumber=0}

\def\majnombres{\ifodd\pageno \else
\pageblanche \fi \pagetitretrue\hautpoly\titreun=0\footnotenumber=0}

\def\hautspages#1#2{\auteurcourant={\ninepcap#1}\titrecourant={\nineit#2}}

\ifnum\chstart=\pageno \pagetitretrue\fi
  

  \def\PAR{\par}
  
  \def\leftnote#1{\vadjust{\setbox1=\vtop{\hsize 20mm\parindent=0pt\eightpoint
  \baselineskip=9pt\rightskip=4mm plus 4mm\vskip-4.7mm#1}\hbox{\kern-2cm\smash{\box1}}}}

  
  \def\raggedcenter{\leftskip=20pt plus 10em  
       \rightskip=\leftskip 
        \parfillskip=0pt 
         \spaceskip=.3333em \xspaceskip=.5em 
          \pretolerance=9999 \tolerance=9999
           \hyphenpenalty=9999 \exhyphenpenalty=9999 }
           
  \def\titrecentre#1{{\parindent0mm\raggedcenter
       \spaceskip=.6em plus .2em minus .2em\xspaceskip=.6em plus .2em minus .2em
        \tit#1\par}}
        


  \def\oui{oui}
  
   \def\fontetitreun{\ifx\paradouze\oui\douzepts\gpdouze\twelvebf\textfont1=\twelveib\else
\quatorzepts\gpquatorze\fourteenbf\fi}

\def\fontetitreunl{\douzepts\textfont1=\twelveib\scriptfont1=\tenib\fourteenti}
 
 \def\fontetitredeux{\textfont1=\eleveni\ifx\paradouze\oui\onzepts\scriptfont1=\ninei\elevenit\else
                        \douzepts\twelveit\fi}
 
   \def\fontetitredeuxb{\ifx\paradouze\oui\onzepts\eleventi\gponze\textfont1=\elevenib\scriptfont1=\nineib
                         \else\douzepts\twelveti\scriptfont1=\twelveib\scriptfont1=\tenib\gpdouze\fi}
                         
\def\fontetitredeuxl{\onzepts\textfont1=\elevenbf\scriptfont1=\ninebf\twelvebf}
  
\def\fontetitretrois{\textfont0=\elevenrm\scriptfont0=\eightrm\textfont1=\eleveni
                      \scriptfont1=\eighti\scriptscriptfont1=\sixi\elevenit}
                      
\def\fontetitrequatre{\textfont0=\elevenrm\scriptfont0=\eightrm\textfont1=\eleveni
                      \scriptfont1=\eighti\scriptscriptfont1=\sixi\elevenrm}
  
  \newcount\titreun\titreun=0
  \newcount\titredeux\titredeux=0
  \newcount\titretrois\titretrois=0
  \newcount\titrequatre\titrequatre=0
  \newcount\enonce\enonce=0
  
  \def\incr#1{\global\advance#1 by 1 {\the #1}}
  \def\avance#1{\global\advance#1 by 1}
  \def\init#1{\global#1=0}
  
  \long\def\Indentation#1#2{\setbox10=\hbox{\fontetitreun#1}
  	                    \ifdim\wd10 < 4mm
                         \setbox10=\hbox to 4mm{\box10\hfill}
                       \else\ifdim\wd10 < 6mm
                         \setbox10=\hbox to 6mm{\box10\hfill}
  	                    \else\ifdim\wd10 < 8mm
                         \setbox10=\hbox to 8mm{\box10\hfill}
                       \else\ifdim\wd10 < 12mm
                         \setbox10=\hbox to 12mm{\box10\hfill}
                       \fi\fi\fi\fi
                       \dimen10=\hsize
                       \advance \dimen10 by -\wd10
                       \noindent \box10 %
                       \ignorespaces
                       \hbox{\vtop{\hsize=\dimen10\raggedright\noindent\fontetitreun#2}}}

  \long\def\paraun#1{\removelastskip\par\medskip\goodbreak\vskip0pt plus.01\vsize\penalty-100
                \vskip0pt plus-.01\vsize
  	              \init{\titredeux}\ifnum\optionparag=1{\init\eqnumber\init\enonce}\else{}\fi
                  \goodbreak{\fontetitreun
  	                \Indentation{\incr{\titreun}.\ }{\fontetitreun #1\par}}\nobreak\medskip}

 %
 %
 \long\def\paraunc#1{\removelastskip\par\bigskip\goodbreak\vskip0pt plus.01\vsize\penalty-100
                \vskip0pt plus-.01\vsize
  	              \init{\titredeux}
                 \ifnum\optionparag=1{\init{\eqnumber}\init\enonce}\else{}\fi
                  \goodbreak
  	                {\parindent0mm\raggedcenter\fontetitreun\incr{\titreun}.\ 
                     \fontetitreun #1\par}\nobreak\medskip}
                     
\newtoks\titreunl
\titreunl={\ifnum\titreun=1{I}\fi%
\ifnum\titreun=2{II}\fi%
\ifnum\titreun=3{III}\fi%
\ifnum\titreun=4{IV}\fi%
\ifnum\titreun=5{V}\fi%
\ifnum\titreun=6{VI}\fi%
\ifnum\titreun=7{VII}\fi%
\ifnum\titreun=8{VIII}\fi%
\ifnum\titreun=9{IX}\fi%
\ifnum\titreun=10{X}\fi%
\ifnum\titreun=11{XI}\fi%
\ifnum\titreun=12{XII}\fi%
\ifnum\titreun=13{XIII}\fi%
}
\long\def\paraunl#1{\removelastskip\par\bigskip\bigskip\goodbreak\vskip0pt plus.01\vsize\penalty-100
                \vskip0pt plus-.01\vsize
  	              \init{\titredeux}\ifnum\optionparag=1{\init\eqnumber\init\enonce}\else{}\fi
                  \goodbreak{\fontetitreunl
  	                \Indentation{\global\advance\titreun by 1{\the\titreunl}.\ }{\fontetitreunl #1\par}}\nobreak\smallskip}

  
  \long\def\paradeux#1{\init{\titretrois}\vskip0pt plus.01\vsize\penalty-10
                \vskip0pt plus-.01\vsize\ifx \elie\oui\medskip\ifnum\titredeux>0\medskip\fi\fi
                 \Indentation{\fontetitredeux\the\titreun${\cdot}$\incr{\titredeux}.}
                              {\fontetitredeux\textfont1=\eleveni#1}\nobreak\par }
  
  \long\def\paradeuxb#1{\init{\titretrois}\vskip0pt plus.001\vsize\penalty-10
                \vskip0pt plus-.01\vsize{\ifx \elie\oui\medskip\ifnum\titredeux>0\medskip\fi\fi
                  \Indentation
  {\fontetitredeuxb\the\titreun${\cdot}$\incr{\titredeux}.}{ \fontetitredeuxb#1}}\nobreak
\smallskip}

\newtoks\titredeuxl
\titredeuxl={\ifnum\titredeux=1{A}\fi%
\ifnum\titredeux=2{B}\fi%
\ifnum\titredeux=3{C}\fi%
\ifnum\titredeux=4{D}\fi%
\ifnum\titredeux=5{E}\fi%
\ifnum\titredeux=6{F}\fi%
\ifnum\titredeux=7{G}\fi%
\ifnum\titredeux=8{H}\fi%
\ifnum\titredeux=9{I}\fi%
\ifnum\titredeux=10{J}\fi%
\ifnum\titredeux=11{K}\fi%
\ifnum\titredeux=12{L}\fi%
\ifnum\titredeux=13{M}\fi%
}
 \long\def\paradeuxl#1{\init{\titretrois}\vskip0pt plus.001\vsize\penalty-10
                \vskip0pt plus-.01
                \vsize \bigskip%
                  \Indentation
     {\fontetitredeuxl\global\advance\titredeux by 1
  \quad \the\titreunl${\cdot}$\the\titredeuxl.}{ \fontetitredeuxl#1}
  \removelastskip\nobreak\smallskip}
  

  \long\def\paratrois#1{\init{\titrequatre}\ifdim\lastskip<\smallskipamount
                \removelastskip\smallskip\fi
                 \vskip0pt plus.01\vsize\penalty-10
                  \vskip0pt
plus-.01\vsize{\ifx \elie\oui\ifnum\titretrois>0\medskip\fi\fi
\Indentation{\fontetitretrois\the\titreun${\cdot}$\the\titredeux${\cdot}$\incr{\titretrois}.\ }
  {\hskip0mm\baselineskip=14pt\fontetitretrois#1}\nobreak\smallskip}}
  
  
  \long\def\paratroisl#1{\init{\titrequatre}\ifdim\lastskip<\smallskipamount
                \removelastskip\fi
                 \vskip0pt plus.01\vsize\penalty-10
                  \vskip0pt
plus-.01\vsize\ifx \elie\oui\bigskip
\fi
\Indentation{\fontetitretrois\quad \quad \the\titreunl{${\cdot}$}\the\titredeuxl${\cdot}$\incr{\titretrois}.\ }
  {\hskip0mm\fontetitretrois#1}\nobreak\smallskip}


  \long\def\paraquatre#1{\ifdim\lastskip<\smallskipamount
                \removelastskip\smallskip\fi
                 \vskip0pt plus.01\vsize\penalty-10
                  \vskip0pt
                  plus-.01\vsize\par
 
\Indentation{\fontetitrequatre \the\titreun{${\cdot}$}\the\titredeux${\cdot}$\the\titretrois${\cdot}$\incr{\titrequatre}.\ }
{\hskip0mm\fontetitrequatre#1}\nobreak\smallskip}


\newtoks\titrequatrel
\titrequatrel={\ifnum\titrequatre=1{a}\fi%
\ifnum\titrequatre=2{b}\fi%
\ifnum\titrequatre=3{c}\fi%
\ifnum\titrequatre=4{d}\fi%
\ifnum\titrequatre=5{e}\fi%
\ifnum\titrequatre=6{f}\fi%
\ifnum\titrequatre=7{g}\fi%
\ifnum\titrequatre=8{h}\fi%
\ifnum\titrequatre=9{i}\fi%
\ifnum\titrequatre=10{j}\fi%
\ifnum\titrequatre=11{k}\fi%
\ifnum\titrequatre=12{l}\fi%
\ifnum\titrequatre=13{m}\fi%
}
\long\def\paraquatrel#1{\ifdim\lastskip<\smallskipamount
                \removelastskip\smallskip\fi
                 \vskip0pt plus.01\vsize\penalty-10
                  \vskip0pt
                  plus-.01\vsize{\bigskip
\Indentation{\global\advance\titrequatre by 1
\fontetitrequatre\quad \quad \quad \the\titreunl${\cdot}$\the\titredeuxl${\cdot}$\the\titretrois${\cdot}$\the\titrequatrel.\ }
{\hskip0mm\fontetitrequatre#1}\nobreak\smallskip}}

\ifx\optionkeys\oui
\def\drefun#1{\definexref{¤#1}{{\the\titreun}}{}} 
\def\drefdeux#1{\definexref{¤#1}{{\the\titreun}.{\the\titredeux}}{}}
\def\dreftrois#1{\definexref{¤#1}{{\the\titreun}.{\the\titredeux}.{\the\titretrois}}{}}
\else
\def\drefun#1{\definexref{prg#1}{{\the\titreun}}{}} 
\def\drefdeux#1{\definexref{prg#1}{{\the\titreun}.{\the\titredeux}}{}}
\def\dreftrois#1{\definexref{prg#1}{{\the\titreun}.{\the\titredeux}.{\the\titretrois}}{}}
\fi

%


  \long\def\propdeux#1#2#3#4{%
       \avance{\enonce}
       \leavevmode\edef\temp{#2}%
         \ifx\temp\empty 
          \else
           \definexref{#2}{#1~{\the\titreun.\the\enonce}}{enonces}
            \definexref{s#2}{{\the\titreun.\the\enonce}}{enonces}
             \fi
\smallskip
      \noindent{\bf#1\ {\bf\the\titreun.\the\enonce{#3}.}\enspace}{\sl#4\par}%
      \ifdim\lastskip<\medskipamount \removelastskip\penalty55\par \fi
   }

  \long\def\propun#1#2#3#4{%
      \avance{\enonce}
       \leavevmode\edef\temp{#2}%
        \ifx\temp\empty 
          \else
           \definexref{#2}{#1~{\the\enonce}}{enonces}
            \definexref{{s#2}}{{\the\enonce}}{enonces}
             \fi
   \par 
     \noindent{\bf#1\ {\bf\the\enonce{#3}.}\enspace}{\sl#4\par}%
     \ifdim\lastskip<\medskipamount \removelastskip\penalty55\medskip\fi
  }
  
  \long\def\prop#1#2#3#4{\ifnum\optionparag=1
                          \propdeux{#1}{#2}{\textfont1=\elevenib#3}{#4} \else\propun{#1}{#2}{\textfont1=\elevenib#3}{#4}\fi}

  \long\def\propt#1#2#3{\ifx\tpf\oui \prop{Th\'eo\-r\`eme}{#1}{#2}{#3}\par
                       \else\prop{Theorem}{#1}{#2}{#3}\par\fi}
  \long\def\Propt#1#2{\propt{#1}{}{#2}}
  \long\def\propl#1#2#3{\ifx\tpf\oui\prop{Lem\-me}{#1}{#2}{#3}\par
                         \else\prop{Lemma}{#1}{#2}{#3}\par\fi}
  \long\def\Propl#1#2{\propl{#1}{}{#2}}
  \long\def\propc#1#2#3{\ifx\tpf\oui\prop{Corol\-laire}{#1}{#2}{#3}\par
                         \else\prop{Corollary}{#1}{#2}{#3}\par\fi}

  \long\def\propd#1#2#3{\ifx\tpf\oui\prop{D\'efi\-nition}{#1}{#2}{#3}\par
                       \else\prop{Definition}{#1}{#2}{#3}\par\fi} 
  
  \long\def\proptd#1#2#3{\ifx\tpf\oui\prop{Th\'eor\`eme et d\'efi\-nition}{#1}{#2}{#3}\par
                       \else\prop{Theorem and definition}{#1}{#2}{#3}\par\fi}


  
  \newcount\optionparag\optionparag=1
  
  \long\def\section#1#2{\ifnum\optionparag=1 \paraun{#2} 
                        \else\goodbreak{\fontetitreun
  	                \Indentation{#1.\ }{#2}}\nobreak\smallskip\fi}

  \def\eqconstruct#1{\ifnum\optionparag=1{\the\titreun\hbox{$\cdot$}#1}\else{#1}\fi}

  
  
  \def\numref{oui}  
  
  \newcount\mesref\mesref=0 
  \def\defbib#1{\ifx\numref\oui\global\advance\mesref by 1 \definexref{#1}{{\the
                 \mesref}}{}\else\definexref{#1}{#1}{}\fi}
  \def\bibtem#1{\defbib{#1}\item{\citer{#1}}}
  \def\citer#1{[\ref{#1}]}
  \def\citeplus#1#2{[\ref{#1}; #2]}

  
  \font\seventeenmsa=msam10 at 17pt    
  \font\fourteenmsa=msam10 at 14pt
  \font\twelvemsa=msam10 at 12pt
  \font\tenmsa=msam10                 
  \font\ninemsa=msam10 at 9pt 
  \font\eightmsa=msam10 at 8pt 
  \font\sevenmsa=msam7 
  \font\sixmsa=msam10 at 6pt
  \font\fivemsa=msam5
  \newfam\msafam\textfont\msafam=\tenmsa\scriptfont\msafam=\sevenmsa\scriptscriptfont\msafam=\fivemsa
  
  \font\seventeenbb=msbm10 at 17pt     
  \font\fourteenbb=msbm10 at 14pt
  \font\twelvebb=msbm10 at 12pt
  \font\tenbb=msbm10                   
  \font\ninebb=msbm10 at 9pt 
  \font\eightbb=msbm10 at 8pt 
  \font\sevenbb=msbm7 
  \font\sixbb=msbm10 at 6pt
  \font\fivebb=msbm5 
  \newfam\bbfam\textfont\bbfam=\tenbb\scriptfont\bbfam=\sevenbb\scriptscriptfont\bbfam=\fivebb
  \def\bb{\fam\bbfam\tenbb}%

  \font\seventeenscaln=eusm10 at 17pt   
  \font\twelvescaln=eusm10 at 12pt
  \font\tenscaln=eusm10                
  \font\ninescaln=eusm10 scaled 900
  \font\eightscaln=eusm10 scaled 800
  \font\sevenscaln=eusm10 scaled 700
  \font\sixscaln=eusm10 scaled 600
   
  \newfam\scalnfam\textfont\scalnfam=\tenscaln\scriptfont\scalnfam=\sevenscaln\scriptscriptfont\scalnfam=\sixscaln
  \def\scaln{\fam\scalnfam\tenscaln}%
  \def\scal{\scaln}
  
  \font\tenscalb=eusb10                

  \font\sevenscalb=eusb10 scaled 700

  \newfam\scalbfam\textfont\scalbfam=\tenscalb\scriptfont\scalbfam=\sevenscalb
  %
  
  %
  %
  \font\fourteenrm=cmr12 scaled 1200
  \font\elevenrm=cmr10 at 11pt
  \font\twelverm=cmr12
  \font\ninerm=cmr9
  \font\eightrm=cmr8      
  \font\sevenrm=cmr7
  \font\sixrm=cmr6

  \font\seventeenpcap=cmcsc10 at 17pt
  \font\tenpcap=cmcsc10                        
  \font\ninepcap=cmcsc9
  \font\eightpcap=cmcsc8
  \font\sevenpcap=cmcsc10 scaled 700
  
  \newfam\pcapfam\textfont\pcapfam=\tenpcap\scriptfont\pcapfam=\sevenpcap
  \def\pcap{\fam\pcapfam\tenpcap}
  
  \font\seventeenrm=cmbx12 scaled 1400

  \font\fourteenbf=cmbx10 scaled 1400
  
  \font\twelvebf=cmbx12
  \font\elevenbf=cmbx10 at 11pt
  \font\ninebf=cmbx9  
  \font\eightbf=cmbx8
  \font\sixbf=cmbx6
  
  \font\tengot=eufm10                           
   
  \font\eightgot=eufm10 at 8truept 
  \font\sevengot=eufm7 
  \font\sixgot=eufm10 at 6 truept 
   
  \newfam\gotfam
  \textfont\gotfam=\tengot\scriptfont\gotfam=\sevengot\scriptscriptfont\gotfam=\sixgot
  \def\got{\fam\gotfam\tengot}%

  
  \def\tit{%
  \textfont0=\seventeenrm\scriptfont0=\tenrm\def\rm{\fam0\seventeenrm}%
  \textfont1=\seventeenib\scriptfont1=\twelveib%
  \textfont2=\seventeensy\scriptfont2=\twelvesy\scriptscriptfont2=\ninesy
  \textfont3=\seventeenex
  \textfont\itfam=\seventeenti
  \def\it{\fam\itfam\seventeenti}%
  \textfont\bbfam=\seventeenbb \scriptfont\bbfam=\twelvebb
  \def\bb{\fam\bbfam\seventeenbb}%
  \textfont\msafam=\seventeenmsa\scriptfont\msafam=\twelvemsa
  \textfont\scalnfam=\seventeenscaln
  \def\pcap{\fam\pcapfam\seventeenpcap}
  \normalbaselineskip=25pt\normalbaselines\rm}

  \font\seventeenti=cmbxti10 scaled 1680
  
  \font\fourteenti=cmbxti10 at 14pt
  
  \font\twelveti=cmbxti10 scaled 1200
  \font\eleventi=cmbxti10 at 11pt

  %
  %
  \font\twelveit=cmti12	
  \font\elevenit=cmti10 scaled 1100
  \font\nineit=cmti9
  \font\eightit=cmti8
  \font\sevenit=cmti7

  %
  %
 
 \font\seventeenib=cmmib10 scaled 1680
  \font\fourteenib=cmmib10 scaled 1400
  \font\twelveib=cmmib10 scaled 1200
  \font\elevenib=cmmib10 scaled 1100
  \font\tenib=cmmib10
\font\eightib=cmmib10 scaled 800
  \font\nineib=cmmib10 scaled 900
\font\sevenib=cmmib10 scaled 700
\font\sixib=cmmib10 scaled 600
\font\fiveib=cmmib10 scaled 500

\ifx\ITAN\oui
\else
\innernewfam\cmmibfam
\textfont\cmmibfam=\tenib
\scriptfont\cmmibfam=\sevenib
\scriptscriptfont\cmmibfam=\fiveib
\def\ib{\fam\cmmibfam\tenib}
\fi

  %
  %
  
  \font\eleveni=cmmi10 scaled 1100
  \font\ninei=cmmi9
  \font\eighti=cmmi8 
  \font\seveni=cmmi7 			                
  \font\sixi=cmmi6
  
  \font\ninesl=cmsl9                    
  \font\eightsl=cmsl8 
  \font\sevensl=cmsl10 at 7pt

  \font\ninett=cmtt9                    
  \font\eighttt=cmtt8
  \font\seventt=cmtt10 scaled 700

  \font\seventeensy=cmsy10 scaled 1680    
  \font\fourteensy=cmsy10 scaled 1400
  \font\twelvesy=cmsy10 scaled 1176
  
  \font\ninesy=cmsy9                      
  \font\eightsy=cmsy8
  \font\sixsy=cmsy6
  \font\seventeenex=cmex10 at 17pt
  \font\fourteenex=cmex10 at 14pt
  \font\twelveex=cmex10 at 12pt
  \font\nineex=cmex10 at 9pt
  \font\eightex=cmex10 at 8pt
  \font\sevenex=cmex10 at 7pt
  \font\sixex=cmex10 at 6pt
  \font\fiveex=cmex10 at 5pt
  
   
  \font\fourteengp=cmmi10 at 14pt
  
  \font\twelvegp=cmmib10 at 12pt
  \font\elevengp=cmmib10 at 11pt
  \font\tengp=cmmib10                          
  \font\ninegp=cmmib10 at 9pt 
  \font\eightgp=cmmib8 
   
  \font\sixgp=cmmib6


  \def\gponze{\textfont0=\elevenbf\scriptfont0=\eightbf\scriptscriptfont0=\sixbf
           \textfont1=\elevengp\scriptfont1=\eightgp\scriptscriptfont1=\sixgp}
  \def\gpdouze{\textfont0=\twelvebf\scriptfont0=\tenbf\scriptscriptfont0=\ninebf
           \textfont1=\twelvegp\scriptfont1=\tengp\scriptscriptfont1=\ninegp}        
  
 \def\gpquatorze{\textfont0=\fourteenbf\scriptfont0=\twelvebf\scriptscriptfont0=\elevenbf
           \textfont1=\fourteengp\scriptfont1=\twelvegp\scriptscriptfont1=\elevengp}

  
  \expandafter\chardef\csname pre amssym.def at\endcsname=\the\catcode`\@
  \catcode`\@=11
  \def\undefine#1{\let#1\undefined}
  \def\newsymbol#1#2#3#4#5{\let\next@\relax
   \ifnum#2=\@ne\let\next@\msafam@\else
   \ifnum#2=\tw@\let\next@\bbfam@\fi\fi
   \mathchardef#1="#3\next@#4#5}
  \def\mathhexbox@#1#2#3{\relax
   \ifmmode\mathpalette{}{\m@th\mathchar"#1#2#3}%
   \else\leavevmode\hbox{$\m@th\mathchar"#1#2#3$}\fi}
  \def\hexnumber@#1{\ifcase#1 0\or 1\or 2\or 3\or 4\or 5\or 6\or 7\or 8\or
   9\or A\or B\or C\or D\or E\or F\fi}
  
  \def\setboxz@h{\setbox\z@\hbox}
  \def\wdz@{\wd\z@}
  \def\boxz@{\box\z@}
  
  \edef\msafam@{\hexnumber@\msafam}
  \mathchardef\dabar@"0\msafam@39
  
  \edef\bbfam@{\hexnumber@\bbfam}
  \def\widehat#1{\setboxz@h{$\m@th#1$}%
   \ifdim\wdz@>\tw@ em\mathaccent"0\bbfam@5B{#1}%
   \else\mathaccent"0362{#1}\fi}
  \def\widetilde#1{\setboxz@h{$\m@th#1$}%
   \ifdim\wdz@>\tw@ em\mathaccent"0\bbfam@5D{#1}%
   \else\mathaccent"0365{#1}\fi}
  \newsymbol\leqq 1335          
  \newsymbol\leqslant 1336
  \newsymbol\lessgtr 1337       
  \newsymbol\backprime 1038     
  \newsymbol\risingdotseq 133A  
  \newsymbol\fallingdotseq 133B 
  \newsymbol\succcurlyeq 133C   
  \newsymbol\geqq 133D          
  \newsymbol\geqslant 133E
  \newsymbol\nmid 232D
  \newsymbol\nexists 2040
  \newsymbol\smallsetminus 2272
  \newsymbol\varnothing 203F
  
  \catcode`\@=\active

  \catcode`\@=11
  \newcount\typofr\typofr=1
  
  \catcode`\;=\active
  \def;{\ifnum\typofr=1\relax\ifhmode\ifdim\lastskip>\z@\unskip\fi
     \kern.2em\fi\string;\else\string;\fi}
  
  \catcode`\:=\active
  \def:{\ifnum\typofr=1\relax\ifhmode\ifdim\lastskip>\z@\unskip\fi
  \penalty\@M\ \fi\string:\else\string:\fi}
  
  \catcode`\!=\active
  \def!{\ifnum\typofr=1\relax\ifhmode\ifdim\lastskip>\z@\unskip\fi
     \kern.2em\fi\string!\else\string!\fi}
  
  \catcode`\?=\active
  \def?{\ifnum\typofr=1\relax\ifhmode\ifdim\lastskip>\z@\unskip\fi
     \kern.2em\fi\string?\else\string?\fi}

  \def\francais{\typofr=1\def\tpf{oui}}
  \def\anglais{\typofr=2\def\tpf{non}\def\english{oui}}
  \def\oui{oui}
  \francais
  
  \catcode`\@=12
  

\ifx\textures\oui
\def\raye #1|{\leavevmode\setbox1=\hbox{#1}%
\raise .5pt\hbox to \wd1{\xleaders\hbox{\rge{$ \sct / $}%
\kern 1pt}\hfill\kern -1pt }\kern -\wd1 \unhbox1\relax }

\def\barre #1|{\leavevmode\setbox1=\hbox{#1}%
\rlap{\Red\vrule height 2.4pt depth -1.2pt width \wd1}\Black \unhbox1\relax}
\else
\def\raye #1|{\leavevmode\setbox1=\hbox{#1}%
\raise .5pt\hbox to \wd1{\xleaders\hbox{\rge{$ \sct / $}%
\kern 1pt}\hfill\kern -1pt }\kern -\wd1 \unhbox1\relax }

\def\barre #1|{\leavevmode\setbox1=\hbox{#1}%
\rlap{\color{red}\vrule height 2.4pt depth -1.2pt width \wd1}\color{black} \unhbox1\relax}

\fi
  

  
  \def\og{\leavevmode\raise.24ex\hbox{$\scriptscriptstyle\langle\!\langle\>$}}    
  \def\fg{\leavevmode\raise.24ex\hbox{$\scriptscriptstyle\>\rangle\!\rangle$}}    

  \def\d{\,{\rm d}}
  \def\dd{{\rm d}}

  \def\z{{\bb Z}}
  \def\r{{\bb R}}
  \def\CC{{\bb C}}

  \def\HH{{\scal H}}

  \def\K{{\scal K}}

  \def\M{{\scal M}}
  
  \def\O{{\scal O}}
  \def\P{{\scaln P}}
  
  \def\R{{\scal R}}

  \def\X{{\scal X}}

  \def\frac#1#2{{#1\over #2}}
  \def\di#1#2{\sct#1\atop{\sct#2}}
  \def\tri#1#2#3{{\sct#1\atop\sct#2}\atop\sct#3}

  \def\numero{n$^{\rm o}\thinspace$}

  \def\qedbox{$\rlap{$\sqcap$}\sqcup$}           
  \def\qed{\nobreak\hfill\penalty250 \hbox{}\nobreak\hfill\qedbox\par }

  \def\numero{n$^{\rm o}\thinspace$}

  \def\plaf#1{\left\lceil#1\right\rceil}
  
  \def\¤{\S\thinspace}

  \def\¥{$\bullet$ }
  
  
  \def\e{{\rm e}}

  \def\epsilon{\varepsilon}

  \def\phi{\varphi}
  \def\theta{\vartheta}
  \def\rho{\varrho}

  \def\dsp{\displaystyle}
  \def\sct{\scriptstyle}
  \def\pf{\noi{\it Proof. }}
  \def\nid{\ifnum\typofr=1\par\noindent{\it D\'emonstration. }\else\pf\fi}
  \def\noi{\noindent}
  \def\rem{\ifnum\typofr=1\noi{\it Remarque.}\ \else\noi{\it Remark.}\ \fi}
  \def\rems{\ifnum\typofr=1\noi{\it Remarques.}\ \else\noi{\it Remarks.}\ \fi}
  \def\re{{\Re e\,}}

  \def\1{{\bf 1}}
  \def\|{\Vert}

  \def\leq{\leqslant}
  \def\geq{\geqslant}

  \def\ie{{i.e.\ }}
  \def\eg{{e.g.}}
  

  \def\fl#1{\left\lfloor #1 \right\rfloor}
  \def\plaf#1{\left\lceil #1 \right\rceil}
  
  \def\sgn{\mathop{\rm sgn}\nolimits}

  \def\log{\mathop{\rm log}\nolimits}



\def\abs#1{\left|#1\right|}


  \def\pmb#1{\setbox0=\hbox{#1}%
  \kern-.025em\copy0\kern-\wd0\kern.05em\copy0\kern-\wd0\kern-.025em\raise .0433em\box0 }

  
  \skewchar\eighti='177 \skewchar\sixi='177
  \skewchar\eightsy='60 \skewchar\sixsy='60
  
  \def\eightpoint{%
  \textfont0=\eightrm\scriptfont0=\sixrm\scriptscriptfont0=\fiverm
  \def\rm{\fam0\eightrm}%
  \textfont1=\eighti\scriptfont1=\sixi
  \scriptscriptfont1=\fivei\def\oldstyle{\fam1\seveni}%
  \textfont2=\eightsy\scriptfont2=\sixsy\scriptscriptfont2=\fivesy
  \textfont3=\eightex\scriptfont3=\sixex
  \textfont\itfam=\eightit
  \def\it{\fam\itfam\eightit}%
  \textfont\slfam=\eightsl
  \def\sl{\fam\slfam\eightsl}%
  \textfont\bbfam=\eightbb \scriptfont\bbfam=\sixbb\scriptscriptfont\bbfam=\fivebb
  \def\bb{\fam\bbfam\eightbb}%
  \textfont\msafam=\eightmsa\scriptfont\msafam=\sixmsa
  \textfont\scalnfam=\eightscaln
  \def\scaln{\fam\scalnfam\eightscaln}
  \textfont\ttfam=\eighttt
  \def\tt{\fam\ttfam\eighttt}%
\textfont\gotfam=\eightgot
  \textfont\bffam=\eightbf\scriptfont\bffam=\sixbf\scriptscriptfont\bffam=\fivebf
  \def\bf{\fam\bffam\eightbf}%
  \ifx\ITAN\oui\else\textfont\cmmibfam=\eightib
       \scriptfont\cmmibfam=\sixib
        \scriptscriptfont\cmmibfam=\fiveib
         \def\ib{\fam\cmmibfam\eightib}
   \fi
  \textfont\pcapfam=\eightpcap
  \def\pcap{\fam\pcapfam\eightpcap}
  \abovedisplayskip=2pt plus2pt minus 2pt
  \belowdisplayskip=2pt plus1pt minus 2pt
  \abovedisplayshortskip= 1pt plus 2pt minus 1pt
  \belowdisplayshortskip= 1pt plus 2pt minus 1pt
  \smallskipamount=2pt plus 1pt minus 2pt
  \medskipamount=3pt plus 2pt minus 2pt
  \bigskipamount=7pt plus 3pt minus 3pt
  \setbox\strutbox=\hbox{\vrule height 5pt depth 2pt width 0pt}%
  \normalbaselineskip=9pt\normalbaselines\rm}

  \def\({\left(}
  \def\){\right)}
  
  \def\footnoterule{\kern -2pt\hrule width 7truecm\kern 2.4pt}
  
  \def\xnotedef#1{\definexref{#1}{\noexpand\number\footnotenumber}{Note}}%

  
  
  \def\ninepoint{%
  \textfont0=\ninerm\scriptfont0=\sixrm\scriptscriptfont0=\fiverm
  \def\rm{\fam0\ninerm}%
  \textfont1=\ninei\scriptfont1=\sixi
  \scriptscriptfont1=\fivei\def\oldstyle{\fam1\ninei}%
  \textfont2=\ninesy\scriptfont2=\sixsy\scriptscriptfont2=\fivesy
  \textfont3=\nineex\scriptfont3=\sixex
  \textfont\itfam=\nineit
  \def\it{\fam\itfam\nineit}%
  \textfont\slfam=\ninesl
  \def\sl{\fam\slfam\ninesl}%
  \textfont\bbfam=\ninebb\scriptfont\bbfam=\sixbb\scriptscriptfont\bbfam=\fivebb
  \def\bb{\fam\bbfam\ninebb}%
  \textfont\msafam=\ninemsa\scriptfont\msafam=\sixmsa\scriptscriptfont\msafam=\fivemsa
  \textfont\scalnfam=\ninescaln
  \def\scaln{\fam\scalnfam\ninescaln}
  \textfont\ttfam=\ninett
  \def\tt{\fam\ttfam\ninett}%
  \textfont\bffam=\ninebf\scriptfont\bffam=\sixbf\scriptscriptfont\bffam=\fivebf
  \def\bf{\fam\bffam\ninebf}%
  \abovedisplayskip=3pt plus2pt minus 2pt
  \belowdisplayskip=3pt plus1pt minus 2pt
  \abovedisplayshortskip= 2pt plus 2pt minus 1pt
  \belowdisplayshortskip= 2pt plus 2pt minus 1pt
  \smallskipamount=2pt plus 1pt minus 2pt
  \medskipamount=3pt plus 2pt minus 2pt
  \bigskipamount=7pt plus 3pt minus 3pt
  \setbox\strutbox=\hbox{\vrule height 5pt depth 2pt width 0pt}%
  \normalbaselineskip=10.5pt plus.3pt minus.3pt\normalbaselines\rm}

  \def\sevenpoint{%
  \textfont0=\sevenrm\scriptfont0=\sixrm\scriptscriptfont0=\fiverm
  \def\rm{\fam0\sevenrm}%
  \textfont1=\seveni\scriptfont1=\sixi
  \scriptscriptfont1=\fivei\def\oldstyle{\fam1\seveni}%
  \textfont2=\sevensy\scriptfont2=\sixsy\scriptscriptfont2=\fivesy
  \textfont3=\sevenex\scriptfont3=\fiveex
  \textfont\itfam=\sevenit
  \def\it{\fam\itfam\sevenit}%
  \textfont\slfam=\sevensl
  \def\sl{\fam\slfam\sevensl}%
  \textfont\bbfam=\sevenbb \scriptfont\bbfam=\sixbb\scriptscriptfont\bbfam=\fivebb
  \def\bb{\fam\bbfam\sevenbb}%
  \textfont\msafam=\sevenmsa\scriptfont\msafam=\sixmsa
  \textfont\scalnfam=\sevenscaln
  \def\scaln{\fam\scalnfam\sevenscaln}
  \textfont\bffam=\sevenbf\scriptfont\bffam=\sixbf\scriptscriptfont\bffam=\fivebf
  \def\bf{\fam\bffam\sevenbf}%
  \textfont\ttfam=\seventt
  \abovedisplayskip=2pt plus2pt minus 2pt
  \belowdisplayskip=2pt plus1pt minus 2pt
  \abovedisplayshortskip= 1pt plus 2pt minus 1pt
  \belowdisplayshortskip= 1pt plus 2pt minus 1pt
  \smallskipamount=2pt plus 1pt minus 2pt
  \medskipamount=3pt plus 2pt minus 2pt
  \bigskipamount=7pt plus 3pt minus 3pt
  \setbox\strutbox=\hbox{\vrule height 5pt depth 2pt width 0pt}%
  \normalbaselineskip=9pt\normalbaselines\rm}

 \def\onzepts{%
 \textfont0=\elevenrm\scriptfont0=\ninerm
 \textfont1=\elevenib\scriptfont1=\ninei}

\def\douzepts{%
  \textfont0=\twelverm\scriptfont0=\tenrm\def\rm{\fam0\twelverm}%
  \textfont1=\twelveib\scriptfont1=\teni%
  \textfont2=\twelvesy\scriptfont2=\tensy\scriptscriptfont2=\eightsy
  \textfont3=\twelveex
  \textfont\itfam=\twelveti
  \def\it{\fam\itfam\twelveti}%
  \textfont\bffam=\twelvebf\scriptfont\bffam=\tenbf\scriptscriptfont\bffam=\eightbf
  \def\bf{\fam\bffam\twelvebf}%
  \textfont\bbfam=\twelvebb \scriptfont\bbfam=\tenbb
  \def\bb{\fam\bbfam\twelvebb}%
  \textfont\msafam=\twelvemsa\scriptfont\msafam=\tenmsa
  \textfont\scalnfam=\twelvescaln
  \normalbaselineskip=15pt\normalbaselines\rm}

\def\quatorzepts{%
  \textfont0=\fourteenrm\scriptfont0=\twelverm\def\rm{\fam0\fourteenrm}%
  \textfont1=\fourteenib\scriptfont1=\twelveib%
  \textfont2=\fourteensy\scriptfont2=\twelvesy\scriptscriptfont2=\tensy
  \textfont3=\fourteenex
  \textfont\itfam=\fourteenti
  \def\it{\fam\itfam\fourteenti}%
  \textfont\bffam=\fourteenbf\scriptfont\bffam=\twelvebf\scriptscriptfont\bffam=\tenbf
  \def\bf{\fam\bffam\fourteenbf}%
  \textfont\bbfam=\fourteenbb \scriptfont\bbfam=\twelvebb
  \def\bb{\fam\bbfam\fourteenbb}%
  \textfont\msafam=\fourteenmsa\scriptfont\msafam=\twelvemsa
  \textfont\scalnfam=\twelvescaln
  \normalbaselineskip=18pt\normalbaselines\rm}


\def\AA{{\it Acta Arith.}}

\def\PLMS{{\it Proc. London Math. Soc.}}

\def\picture #1 by #2 (#3){\leavevmode\vbox to #2{
     \hrule width #1 height 0pt depth 0pt
      \vfill
       \special{picture #3}}}

\def\illustration #1 by #2 (#3) scaled #4{\dimen1=#2
  \divide\dimen1 by 1000
  \multiply\dimen1 by #4
  \vtop to \dimen1{\dimen1=#1
  \divide\dimen1 by 1000
  \multiply\dimen1 by #4
  \hsize=\dimen1\vss
  \noindent\special{illustration #3 scaled #4}}}

 \fi
\ifx\optionkeymacros\undefined\else \fi

\catcode`\Œ=\active\defŒ{{\aa}}       
\catcode`\º=\active\defº{\int}        
\catcode`\=\active\def{\c c}        
\catcode`\¶=\active\def¶{\partial}    
\catcode`\Ä=\active\defÄ{\oint}       
\catcode`\Æ=\active\defÆ{\triangle}   
\catcode`\Â=\active\defÂ{\neg}        
\catcode`\µ=\active\defµ{\mu}         
\catcode`\¿=\active\def¿{{\o}}        
\catcode`\¹=\active\def¹{\pi}         
\catcode`\Ï=\active\defÏ{{\oe}}       
\catcode`\§=\active\def§{{\ss}}       
\catcode`\ =\active\def {\dagger}     
\catcode`\Ã=\active\defÃ{\sqrt}       
\catcode`\·=\active\def·{\Sigma}      
\catcode`\Å=\active\defÅ{\approx}     
\catcode`\½=\active\def½{\Omega}      
\catcode`\£=\active\def£{{\it\$}}     
\catcode`\°=\active\def°{\infty}      
\catcode`\¤=\active\def¤{{\S}}        
\catcode`\¦=\active\def¦{{\P}}        
\catcode`\¥=\active\def¥{\bullet}     
\catcode`\»=\active\def»{\leavevmode\raise.585ex\hbox{\b a}}      
\catcode`\¼=\active\def¼{\leavevmode\raise.6ex\hbox{\b o}}        
\catcode`\­=\active\def­{\not=}       
\catcode`\²=\active\def²{\leq}        
\catcode`\³=\active\def³{\geq}        
\catcode`\Ö=\active\defÖ{\div}        
\catcode`\É=\active\defÉ{{\dots}}     
\catcode`\¾=\active\def¾{{\ae}}       
\catcode`\Ç=\active\defÇ{\og}         
\catcode`\Ò=\active\defÒ{``}          
\catcode`\Á=\active\defÁ{!`}          
\catcode`\¢=\active\def¢{\rlap/c}     
\catcode`\Ô=\active\defÔ{`}           
\catcode`\Õ=\active\defÕ{'}           


\catcode`\=\active\def{{\AA}}       
\catcode`\'=\active\def'{\c C}        
\catcode`\¯=\active\def¯{{\O}}        
\catcode`\¸=\active\def¸{\Pi}         
\catcode`\Î=\active\defÎ{{\OE}}       
\catcode`\®=\active\def®{{\AE}}       
\catcode`\×=\active\def×{\diamond}    
\catcode`\¡=\active\def¡{\accent'27}  
\catcode`\Ó=\active\defÓ{''}          
\catcode`\±=\active\def±{\pm}         
\catcode`\È=\active\defÈ{\fg}         
\catcode`\À=\active\defÀ{?`}          
\catcode`\Ð=\active\defÐ{--}          
\catcode`\Ñ=\active\defÑ{---}         


\catcode`\Š=\active\defŠ{\"a}        
\catcode`\'=\active\def'{\"e}        
\catcode`\•=\active\def•{\"{\i}}     
\catcode`\š=\active\defš{\"o}        
\catcode`\Ÿ=\active\defŸ{\"u}        
\catcode`\Ø=\active\defØ{\"y}        
\catcode`\å=\active\defå{\^A}        
\catcode`\€=\active\def€{\"A}        
\catcode`\…=\active\def…{\"O}        
\catcode`\†=\active\def†{\"U}        
\catcode`\‡=\active\def‡{\'a}        
\catcode`\Ž=\active\defŽ{\'e}        
\catcode`\'=\active\def'{\'{\i}}     
\catcode`\—=\active\def—{\'o}        
\catcode`\œ=\active\defœ{\'u}        
\catcode`\ƒ=\active\defƒ{\'E}        
\catcode`\æ=\active\defæ{\^E}        
\catcode`\é=\active\defé{\`E}        %
\catcode`\ˆ=\active\defˆ{\`a}        
\catcode`\=\active\def{\`e}        
\catcode`\"=\active\def"{\`{\i}}     
\catcode`\˜=\active\def˜{\`o}        
\catcode`\=\active\def{\`u}        
\catcode`\Ë=\active\defË{\`A}        
\catcode`\‹=\active\def‹{\~a}        
\catcode`\–=\active\def–{\~n}        
\catcode`\›=\active\def›{\~o}        
\catcode`\Ì=\active\defÌ{\~A}        
\catcode`\"=\active\def"{\~N}        
\catcode`\Í=\active\defÍ{\~O}        
\catcode`\‰=\active\def‰{\^a}        
\catcode`\=\active\def{\^e}        
\catcode`\"=\active\def"{\^{\i}}     
\catcode`\™=\active\def™{\^o}        
\catcode`\ž=\active\defž{\^u}        

\let\optionkeymacros\null

\def\paradouze{oui}

\optionparag=1
\def\paradouze{oui}
\anglais

    \font\tenrsfs=rsfs7 at 10pt

    \font\sevenrsfs=rsfs7
    \font\sixrsfs=rsfs7 at 6pt
    
\newfam\rsfsfam\textfont\rsfsfam=\tenrsfs\scriptfont\rsfsfam=\sevenrsfs\scriptscriptfont\rsfsfam=\sixrsfs
    \def\rsfs{\fam\rsfsfam\tenrsfs}%
\def\pnu{p^\nu}

\def\gB{{\got B}}

\def\gf{{\got f}}

\def\gtg{{\got g}}

\def\FF{{\scal F}}
\def\GG{{\rsfs G}}

\def\HH{{\scal H}}
\def\Sc{{\scal S}}
\def\X{{\scal X}}
 \def\pnu{{p^\nu}}
\def\fl#1{\left\lfloor #1 \right\rfloor}

\ifx\montrerlabels\oui

\input montrerlabels.tex
\fi

\dimstand
\vsize235truemm
\hautspages{R. de la Bretche \& G. Tenenbaum}{Remarks on the Selberg--Delange method}
\dateheure

\titrecentre{Remarks on the Selberg--Delange method}
\bigskip\medskip
\centerline{RŽgis de la Bretche \& GŽrald Tenenbaum}
\bigskip\bigskip
{\eightpoint\leftskip1cm\rightskip1cm
\noi{\bf Abstract.} Let $\varrho$ be a complex number and let $f$ be a multiplicative arithmetic function whose Dirichlet series takes the form $\zeta(s)^\varrho G(s)$, where $G$ is associated to a multiplicative function~$g$. The classical Selberg-Delange method furnishes asymptotic estimates for the averages of $f$ under assumptions of either analytic continuation for $G$, or absolute convergence of a finite number of derivatives of $G(s)$ at $s=1$. We consider different set of hypotheses, not directly comparable to the previous ones, and investigate how they can yield  sharp  asymptotic estimates for the averages of~$f$. \PAR
\medskip\noi
{\bf Keywords: \rm Averages of multiplicative functions, Selberg-Delange method, Dirichlet series, powers of the Riemann zeta function.}  \par
\smallskip 
\noi \bf 2020 Mathematics Subject Classification: \rm  11N37.\par }
\par 
\bigskip\bigskip
\paraun{Introduction and statement of results}
In a series of papers published in 1953 and 1954, \citer{Sa53}, \citer{Sa54}, L. G. Sathe studied the local laws of the arithmetic functions counting the number of prime factors,  with or without multiplicity, a problem previously considered by Hardy and Ramanujan.  Sathe's results provided asymptotic formulae while only upper and lower bounds were previously known.  However Sathe's method, based on induction formulae, involved very long and technical estimates. In the same year 1954, Selberg devised a fruitful method based on the idea that the Dirichlet series of the search for probabilities may be expressed through   Taylor coefficients of powers of the Riemann zeta function. This idea was then systematically developed by Delange \citer{De59}, \citer{De71}. In the second author's book \citer{Te15} (latest edition, first  in 1990), the results were generalized and made effective regarding various parameters---the method being there named after Selberg and Delange.
\par 
 This theory provides estimates for counting functions associated to  Dirichlet series of the form $$F(s)=\zeta(s)^\varrho G(s),\eqdef{hypSD}$$ where $\varrho$ is a complex number, $\zeta(s)$ is Riemann's zeta function, and $G(s)$ satisfies suitable regularity conditions. Usually, $G$ is associated with a multiplicative arithmetic function $g$, but this need not be so --- see, \eg, \citer{HTW08}. In the sequel, we shall however concentrate on the case when $g$ is indeed multiplicative. Then $F(s)$ is the Dirichlet series of a multiplicative function~$f$.\par 
In  \citer{Te15}, two types of assumptions on $G(s)$ are considered : (a) analytic continuation at the left of the vertical line $\sigma=\re s=1$; (b) absolute convergence at $s=1$ for a finite number of derivatives.
\par 
In a recent work \citer{GK19}, Granville and Koukoulopoulos, propose a third type of condition, implying mere convergence, instead of absolute convergence, for a finite number of right-derivatives of $G(s)$ at $s=1$ . However, their analysis actually rests upon  a much stronger assumption, viz., for some constant $A>0$, not necessarily an integer,
$$\sum_{p\leqslant x}g(p)\log p=\sum_{n\leqslant x}\{f(p)-\varrho\}\log p\ll x/(\log x)^A\quad(x\geqslant 2).\eqdef{hypg}$$
 Some extra, secondary hypotheses are also needed for the values $f(\pnu)$ at prime powers.  Here and in the sequel, the letter $p$ denotes a prime number. For the sake of comparison, it is worthwhile to note that hypothesis~(b) above essentially amounts to $\sum_{p}|g(p)|(\log p)^j/p<\infty$ for a finite number of exponents $j$.
\par 
Similar, but weaker, conditions  have been considered by Wirsing \citer{Wi61}, in the frame of comparison theorems,  evaluating the ratio of averages of $f$ and of a non-negative majorant. Via a further weakening of the hypotheses, such results were improved in \citer{Wi67}. In \citer{Te17}, the second author considered generalizations, and obtained effective forms of the results. In the present work, the viewpoint is  quite different: taking advantage of the strength of assumptions like \eqref{hypg}, one aims at directly deriving  an asymptotic estimate for the averages of $f$. As in the classical Selberg-Delange approach, this also enters in the frame of comparison theorems, but the average of $f(n)$ is now compared with that of $\tau_\varrho(n)$---the $n$th coefficient in the Dirichlet series expansion of~$\zeta(s)^\varrho$--- instead of being compared with that of a majorant.
\par \goodbreak
 Assuming \eqref{hypSD}, \eqref{hypg}, and $|f|\leqslant \tau_r$ for some parameter $r>0$,  the main result in \citer{GK19} states that, with $J:=\plaf{A-1}$ and  suitable coefficients $\{\lambda_j(f)\}_{j=0}^{J}$, we have, for $x\geqslant 3$, 
$$M(x;f):=\sum_{n\leqslant x}f(n)=x(\log x)^{\varrho-1}\bigg\{\sum_{0\leqslant j\leqslant J}{\lambda_j(f)\over (\log x)^j}\bigg\}+O\bigg({x(\log_2x)^{\delta_{A,J+1}}\over (\log x)^{A+1-r}}\bigg),\eqdef{estGK}$$
with Kronecker's $\delta$-notation. Here and henceforth, $\log_k$ denotes the $k$-fold iterated logarithm.\par 
\goodbreak
Formula \eqref{estGK} is specially interesting when $A$ is small, for less is then required on $g$. However, it furnishes no more than an upper bound when $A\leqslant r-\re\varrho$. Incidentally, under the weaker assumption that the series
$$\sum_{p}{g(p)\over p}$$ converges, and arguing as in \citer{HT91} (see also \citeplus{Te15}{th. III.4.14}),
theorem 1.1 of \citer{Te17} readily yields, for real $f$,
$$M(x;f)\ll x(\log x)^{r-1-\min\{1,K(r-\varrho)\}}\qquad (x\geqslant 2),\eqdef{majHT}$$
where $K\approx 0.32867$ is optimal. Moreover, in the case $\varrho=0$, the same technique furnishes
$$M(x;f)\ll x(\log x)^{r-1-\min\{1,(1-2/\pi)r\}}\qquad (x\geqslant 2).\eqdef{majHT+}$$
In particular, \eqref{majHT} supersedes \eqref{estGK} as soon as $A<K(r-\varrho)\leqslant 1$. Analogous upper bounds are available for complex $f$, under suitable hypotheses upon $f(p)$: see \eg\ \citer{Ha95}.
\par \medskip
The purpose of the present work is two-fold: (a) to investigate refinements of \eqref{hypg} enabling an improvement of the error term  in \eqref{estGK} by replacing the exponent $r$ of $\log x$ by  $\re\varrho$, as expected in view of standard estimates in the theory; (b) to propose a simpler and more natural (i.e. relying on a direct convolution argument) proof of \eqref{estGK}, with possibly weaker hypotheses. 
\par 
 In the first direction, we meet the set goal  when $\re\varrho\geqslant 0$,  where $\varrho$ is the exponent appearing in the generic assumption \eqref{hypSD}. This  latter restriction is actually necessary for achieving target~(a):  we construct  a family of counter-examples in Section \refn{prgcontrex1}. \par 
  In the following statement, hypothesis \eqref{hypg} is replaced by a short interval version, with the same value of~$A$. The other assumptions concern $|f|$: we use those now classical introduced by Shiu \citer{Sh80}, although they  could be somewhat weakened if needed. Accordingly, we  define  the class $\Sc(B)$ of those multiplicative functions $f$ such that:
\par \medskip
{\parindent 10mm
(i) \qquad $|f(\pnu)|\leqslant B^\nu$\quad$(\nu\geqslant 0)$,\par \smallskip
(ii) \qquad $\forall \varepsilon>0\ \exists C=C_\varepsilon: |f(n)|\leqslant C_\varepsilon n^\varepsilon$ $(n\geqslant 1)$,}
\medskip
\noi and, for $r>0$, we consider the subclass $\Sc(B,r)$ of those $f$ satisfying the extra assumption:\medskip
{\parindent 10mm
(iii) \qquad $\dsp \sum_{p\leqslant x}{|f(p)|\over p}\leqslant r\log_2x+O(1)\quad  (x\geqslant 3).$}
\Propt{SD+}{Let $A>0$, $B>0$, $0<\alpha<1$, $r>0$, $\varrho\in\CC$, $J:=\plaf{A-1}$, and let $f\in\Sc(B,r)$ verify 
$$\sum_{x<p\leqslant x+z}f(p)\log p=\varrho z+O\Big({z\over (\log x)^A}\Big) \quad(x\geqslant 2,\,   x^{1-\alpha}\leqslant z\leqslant x). \eqdef{hyp}$$
 Then, for suitable constants $\{\lambda_j(f) \}_{0\leqslant j\leqslant J}$ and with $\vartheta:=B+|\varrho|+ 1$, we have
$$M(x;f)=x(\log x)^{\varrho-1}\sum_{0\leqslant j\leqslant J}{ \lambda_j(f) \over (\log x)^j}+O\bigg({x(\log_2x)^{\vartheta}\over (\log x)^{A+1-\max(0,\re\varrho)}}\bigg)\qquad (x\geqslant 3).\eqdef{fSD}$$
}
\smallskip
\goodbreak
The coefficients $ \lambda_j(f) $ may be described as follows. Representing $f=\tau_\varrho*g$  consistently with~\eqref{hypSD}, we prove the convergence of the series
$$ \gamma_j(g) :=\sum_{n\geqslant 1}{g(n)(\log n)^j\over n}\qquad (0\leqslant j\leqslant J),\eqdef{defgamj}$$
and derive
$$\lambda_j(f)={1\over \Gamma(\varrho-j)}\sum_{\ell+h=j}{\alpha_\ell(\varrho) {\gamma_h(g)}\over \ell!h!}\qquad (0\leqslant j\leqslant J),\eqdef{flj}$$
where $\Gamma$ is Euler's function and $\alpha_\ell(\varrho)/\ell!$ is the $\ell$-th Taylor coefficient at the origin of $$\{s\zeta(s+1)\}^\varrho/(s+1).$$ Alternatively, we also have
$$\lambda_j(f)={\dd^j \{s^\rho F(s+1)/(s+1)\}\over j!\Gamma(\varrho-j)\d s^j}(0+)\qquad (0\leqslant j\leqslant J).\eqdef{flambdaj}$$
In particular, $\lambda_j(f)=0$ for all $j$ if $\varrho$ in an integer $\leqslant 0$.
\par
\goodbreak
As for our second aim, namely goal (b) described above, we apply {\it friable convergence} (see below) to show the following result,  essentially equivalent, in the intersection of the respective validity domains, to \hbox{\citeplus{GK19}{th. 1}}.  For  $r>0$, $\sigma\in]0,1[$, we define the class $\FF(r,\sigma)$ comprising those complex multiplicative functions $f$ such that
$$\eqalign{&\sum_{v<p\leqslant w}{|f(p)|\over p}\leqslant r\log \Big({\log w\over \log v}\Big)+O(1)\quad(w\geqslant v\geqslant 2),\cr
&\sum_{p}\Big\{{|f(p)|^2\over p^{2\sigma}}+\sum_{\nu\geqslant 2}{|f(\pnu)|\over p^{\nu\sigma}}\Big\}<\infty.\cr}\eqdef{FFB}$$
These conditions are weaker than those of \citeplus{GK19}{th. 1} and not directly comparable to those described in \citeplus{GK19}{\S\thinspace7}.  For instance, letting $p_k$ denote the $k$th prime number, conditions \eqref{FFB} allow $f(p_{k^2})\asymp k/\log 2k$, whereas the latter do not.\par 
Contrary to that of \citer{GK19}, our analysis does not ascribe a special role to the case when $A$ is an integer. 
\Propt{SDf}{Let $A>0$, $J:=\plaf{A-1}$, $r>0$, $\varrho\in\CC$, $\sigma\in]0,1[$, and let $f\in\FF(r,\sigma)$ be a multiplicative function such that  $p\mapsto g(p):=f(p)-\varrho$ verifies \eqref{hypg}.
 Then, for suitable constants  $\beta>0$ and  $\{\lambda_j(f)\}_{0\leqslant j\leqslant J}$, we have
$$M(x;f)=x(\log x)^{\varrho-1}\sum_{0\leqslant j\leqslant J}{\lambda_j(f)\over (\log x)^j}+O\bigg({ x(\log_3x )^{\beta}\over (\log x)^{A+1-r}}\bigg).\eqdef{fSDf}$$
}
We shall see that the $\lambda_j(f)$ are still given by \eqref{flambdaj}. \par \goodbreak
 In Section \refn{prgcontrex2}, we show that, even when $\re\varrho\geqslant 0$, the remainder term of \eqref{fSDf} cannot be sharpened so as to meet that of \eqref{fSD} ---although a possibility of some improvement  remains open if $f$ is real. To this end, we construct a family of counterexamples $f\in\FF(\sigma,r)$ satisfying \eqref{hypg} with $\varrho=0$,  $0<A<r$ (resp. $0<A<2r/\pi$ in the real case), but contravening the short interval condition \eqref{hyp}, and for which the exponent $r$ in \eqref{fSDf} cannot be replaced by $cr$ if $c<1$ (resp. $c<2/\pi$).   \par 
 As a final remark, we note the following: as that of \citer{GK19}, our analysis heavily depends on the asymptotic expansion for
 $$T_\varrho(x):=M(x;\tau_\varrho)=\sum_{n\leqslant x}\tau_\varrho(n)\eqdef{defTrho}$$
  provided by the Selberg-Delange method;  since the arithmetical functions under consideration are of the form $f=g*\tau_\varrho$, refinements regarding hypotheses on $g$ or on its generating series $G(s)$ may hence be regarded as genuine parts of this theory. 
\bigskip
\paraun{On the case $\varrho=r$}
The bound \eqref{majHT} shows that the error term of \eqref{estGK} is not optimal in the general case. Moreover, \ref{SD+}  provides fairly standard assumptions under which an essentially optimal remainder is achieved. However, inasmuch the power of the logarithm is concerned,~\eqref{estGK} is expected to be sharp when $\varrho=r$. This latter case is discussed  in \citer{GK19}, where it is  confirmed that the error term   may be replaced by a quantity  $\asymp x(\log x)^{r-A-1}$ when $\varrho=r\geqslant 1$, $A\notin\z$. 
\par 
Assuming $\varrho=r$ essentially amounts to considering $f\geqslant 0$.  In this latter circumstance, an asymptotic {\it formula} with optimal remainder may be obtained in a very simple way under lighter hypotheses. We present the details below.
\Propt{thpos}{Let $A>0$, $\sigma\in]0,1[$, $r>0$, and assume  $f$ is a non-negative multiplicative function such that \medskip 
{\parindent 10mm
\rm(i)\qquad $\dsp\sum_{p\leqslant x}f(p)\log p=rx+O\Big({x\over (\log x)^A}\Big)$ $\ (x\geqslant 2)$,
\par 
\rm(ii)\qquad  $\dsp\sum_{p}\bigg\{{f(p)^2\over p^{2\sigma}}+\sum_{\nu\geqslant 2}{f(\pnu)\over p^{\nu\sigma}}\bigg\}<\infty$.\par }
\smallskip
We then have
$$M(x;f)=\lambda_0(f)x(\log x)^{r-1}\Big\{1+O\Big({(\log_2x)^{\delta_{1,A}}\over (\log x)^{\min(1,A)}}\Big)\Big\}\quad(x\geqslant 2),\eqdef{fatr}$$
with $\dsp\lambda_0(f) :={1\over \Gamma(r)}\prod_{p}\Big(1-{1\over p}\Big)^r\sum_{\nu\geqslant 0}{f(\pnu)\over \pnu}\cdot$}
\nid  We have
$$\eqalign{M(x;f)\log x&=\sum_{n\leqslant x}f(n)\log n+\sum_{n\leqslant x}f(n)\log (x/n)\cr
&=\sum_{m\leqslant x}f(m)\sum_{\di{\pnu\leqslant x/m}{p\,\nmid\, m}}f(\pnu)\log \pnu+\int_1^x {M(t; f)\over t}\d t\cr
&=r x\sum_{m\leqslant x}{f(m)\over m}+O\bigg(R+S+x(\log x)^{r-1}\bigg),\cr}\eqdef{convol}$$
with
$$\eqalign{R&:=\sum_{\di{mp\leqslant x}{p\,\mid\,m}}f(m)f(p)\log p+\sum_{\di{\pnu m\leqslant x}{\nu\geqslant 2}}f(m)f(\pnu)\log \pnu,\quad
S:=\sum_{m\leqslant x}{xf(m)\over m(\log 2x/m)^A},\cr}$$
and where the last integral has been estimated by an  Halberstam-Richert type bound---see, \eg, \citeplus{Te15}{th.~{III}.3.5}.\par 
\goodbreak
Now, with $\alpha:=(\sigma+1)/2 <1 $, we have
$$\eqalign{R&\leqslant \sum_{m\leqslant x}f(m)\sum_{\di{\pnu m\leqslant x}{\nu\geqslant 2}}\{f(p)f(p^{\nu-1})+f(\pnu)\}\log \pnu\cr
&\ll
\sum_{m\leqslant x}{f(m)x^\alpha\over m^\alpha}\sum_{p}\sum_{\nu\geqslant 2}{f(p)f(p^{\nu-1})+f(\pnu)\over p^{\nu\sigma}}\ll
\sum_{m\leqslant x}{f(m)x^\alpha\over m^\alpha}\ll x(\log x)^{r-1},\cr}$$
where we used assumption (ii). Moreover, partial summation furnishes
  $$S\ll x(\log x)^{r-1-\min(1,A)}(\log_2x)^{\delta_{1,A}}.$$
It remains to note that, by \citeplus{TW08}{th. 3.3} (with $\kappa=r$ and $y=x$), we have
$$\sum_{m\leqslant x}{f(m)\over m}=\Big\{1+O\Big({1\over \log x}\Big)\Big\} {\lambda_0(f)\over r}(\log x)^{r}\qquad(x\geqslant 2),$$
 and carry back into \eqref{convol}. 
\vskip-5mm\qed
\goodbreak\medskip
\bigskip
\paraun{Proof of \ref{SD+}}
\paradeuxb{Preparation}
Given  $\varrho\in\CC$, let $\M(\varrho)$ denote the class of those multiplicative functions $f$ satisfying
$$f(\pnu)=\tau_\varrho(\pnu)+\tau_\varrho(p^{\nu-1})\{f(p)-\varrho\}\qquad (p\geqslant 2,\,\nu\geqslant 1).\eqdef{red}$$
 When $f\in\Sc(B,r)$, we can write $f=\gf*h$ where $\gf\in\M(\varrho)$, $h$ is multiplicative, supported on the set of squareful integers, and such that the series $$\sum_{n\geqslant 1}{h(n)(\log n)^j\over n}\qquad (j\geqslant 0)$$
are absolutely convergent and uniformly bounded in terms of $B$ for  bounded $j$. This is enough to yield the required estimate, arguing as in \citeplus{Te15}{th. II.5.4}. Therefore, it will be sufficient to prove that \eqref{fSD} holds  when $f\in\M(\varrho)$, $|f(p)|\leqslant B$, and $f$ satisfies \eqref{hyp}. 
   \par 
As a consequence of this reduction, we can assume that $f=\tau_\varrho*g$, where $g$ is supported on the set of squarefree integers. Moreover, $g\in\Sc(\gB,r+|\varrho|)$ with $\gB:=B+|\varrho|$, and satisfies~\eqref{hypg}. 
\medskip\goodbreak
\paradeuxb{Main estimates}
The main part of the argument consists in providing effective estimates for averages of the function $g$ defined above.  
\Propl{Gjk}{For a suitable constant $c_0$ and uniformly for $x\geqslant 3,\,0\leqslant j\leqslant J,\,k\ll\log_2x$, we have
$$\eqalign{G_{j,k}(x)&:=\sum_{\di{n\leqslant x}{\omega(n)=k}}g(n)(\log n)^j\ll{x(\gB\log_3x+c_0)^{k-1}\over (k-1)!(\log x)^{A+1-j}}\cdot\cr}\eqdef{majgamjk}$$
}
\par 
For further reference we note right away that a consequence of \eqref{majgamjk} is that 
$$\leqalignno{G_j(x)&:=\sum_{n\leqslant x}g(n)(\log n)^j\ll_j{x(\log_2x)^{\gB}\over (\log x)^{A+1-j}}\qquad (j\geqslant 0,\,x\geqslant 3),&\eqdef{majGj}\cr
\gtg_j(x)&:=\sum_{n\leqslant x}{g(n)(\log n)^j\over n}=\gamma_j(g)+O\bigg({(\log_2x)^{\gB}\over (\log x)^{A-j}}\bigg)\quad(0\leqslant j\leqslant J,\,x\geqslant 3).&\eqdef{omsomgj}\cr}$$
\par 
Indeed, it suffices to apply, \eg, lemma 1 from \citer{Te00} in order to note that, for large, constant~$D$, the contribution of those $n$ with \hbox{$\omega(n)>D\log_2x$} is negligible and  then appeal to \eqref{majgamjk} for  $k\leqslant D\log_2x$.
\par 
 In the following proof and henceforth, we let $P^+(n)$---resp. $P^-(n)$---denote the largest---resp. the smallest---prime factor of an integer $n>1$, and make the standard convention that $P^+(1)=1$,\ $P^-(1)=\infty$. 
\medskip\goodbreak  
\noi{\it Proof of \ref{Gjk}.} By partial summation, it is enough to consider $j=0$. We may also plainly restrict to bounding the subsum over $n\in ]\sqrt{x},x]$.
\par 
Given a suitably large parameter $\K$,  set $\X:=(\log x)^\K$ and represent each $n$ arising in \eqref{majgamjk} as $n=md$ with $P^+(m)\leqslant \X$, $P^-(d)>\X$. By, for instance, \citeplus{Te87}{lemma 2} (a Rankin-type bound), we see that the contribution  to $G_{0,k}(x)$  of $d\leqslant x^{1/4}$ is, for suitable positive constants~$c_j$, 
$$\ll\sum_{d\leqslant x^{1/4}}{|g(d)|}\sum_{\di{x^{1/4}<m\leqslant x/d}{P^+(m)\leqslant \X}}{|g(m)|}\ll\sum_{d\leqslant x^{1/4}}{|g(d)|\over d}x^{1-c_1/\log \X}\ll x^{1-c_2/\log \X}.$$
Whence
$$\eqalign{
G_{0,k}(x)=\sum_{\di{s+t=k}{t\geqslant 1}}\sum_{\tri{m\leqslant x^{3/4}}{\omega(m)=s}{P^+(m)\leqslant \X}}g(m)\GG_t\Big({x\over m};\X\Big)+O\Big(x^{1-c_3/\log_2x}\Big),\cr}\eqdef{GGh1}$$
with
$$\GG_t(w;\X):=\sum_{\tri{x^{1/4}<d\leqslant w}{\omega(d)=t}{P^-(d)>\X}}g(d)\quad\big(1\leqslant t\leqslant k,\,x^{1/4}<w\leqslant x\big).\eqdef{defHHt}$$
Let $\HH:=\alpha \K$, where $\alpha>0$ is the constant appearing in \eqref{hyp}, and let $\delta:=1/(\log x)^\HH$. Put $I(\ell):=]\e^{\delta\ell},\e^{\delta(\ell+1)}]$ $(\ell>L:=(\K\log_2x)/\delta)$. The contribution to \eqref{defHHt} from those integers $d$ having at least two prime factors in a same interval $I(\ell)$ is
$$\ll \gB^t\sum_{1\leqslant j<t}\sum_{\tri{\X<p_1<\cdots<p_{t}}{p_1\cdots p_t\leqslant w}{p_{j+1}<p_j\e^\delta}}{w\over p_1\cdots p_t}\ll \delta w(\log w)^{\gB}\ll{w\over (\log x)^{A+2}},$$
for a suitable choice of $\K$ and hence of $\HH$. 
\par \goodbreak
For the remaining integers, split the prime factors of the summation variable $d$ in \eqref{defHHt} among the various $I(\ell)$ and consider the multiple sum over all  hypercubes that are hit. 
If $\ell_1< \ell_2< \ldots<  \ell_t$ is a sequence of admissible indexes, then
$\delta(\ell_1+\ldots+\ell_t)\leqslant \log w$ since $d\leqslant w$ in~\eqref{defHHt}.
If some product $v:=\prod_{j=1}^tp_j$ happens to be $> w$, we must have
$w<v\leqslant  w\e^{t\delta}.$
Therefore the total contribution of those admissible hypercubes containing at least one product $> w$~is
$$\ll\sum_{{w<v\leqslant  w\e^{t\delta}}}{\gB^{\omega(v)}\mu(v)^2}\ll w(\log w)^{\gB-1}\delta t\ll{w\over (\log x)^{A+2}}\cdot$$
\par 
We have shown so far that, for relevant values of $t$ and $w$, 
$$\GG_t(w;\X)=S_t+O\bigg({w\over (\log x)^{A+2}}\bigg),$$ 
with, for some absolute constants $C_j$ $(j=0,1)$,
$$\eqalign{S_t&:=\sum_{\di{L<\ell_1<\cdots<\ell_t}{(\log x)/4\delta<\ell_1+\ldots+\ell_t\leqslant (\log w)/\delta}}\prod_{1\leqslant j\leqslant t}\sum_{p\in I(\ell_j)}g(p)\cr&\ll\sum_{\di{L<\ell_1<\cdots<\ell_t}{(\log x)/4\delta<\ell_1+\ldots+\ell_t\leqslant (\log w)/\delta}}\prod_{1\leqslant j\leqslant t}{C_0\e^{\ell_j\delta}\over \delta^A\ell_j^{A+1}}\cr
&\ll{C_1^tt^{A+1}\over  (t-1)! (L\delta)^{A(t-1)}}\sum_{(\log x)/4\delta<\ell\leqslant (\log w)/\delta}{\e^{\delta\ell}\over \delta^A\ell^{A+1}}\cr
&\ll{C_1^tt^{A+1}w\over (t-1)!(\delta L)^{A(t-1)} (\log x)^{A+1}}\ll{w\over (t-1)!(\kappa\log_2x)^{t-1}(\log x)^{A+1}},\cr}$$
where $\kappa=\kappa(\K)$ may be chosen as large as we wish. Here, we made use of the short-interval  assumption \eqref{hyp}  at the very first step.  In the third line, 
we set  $\ell:= \ell_1+\ldots+\ell_t$ and bounded   $1/\ell_t^{A+1}$ by $t^{A+1}/\ell^{A+1}$.  
Carrying back into \eqref{GGh1}, we obtain
$$G_{0,k}(x)\ll\sum_{\di{s+t=k}{t\geqslant 1}}{x(\gB\log_3x+c_3)^s\over s!(t-1)!(\kappa\log_2x)^{t-1}(\log x)^{A+1}}\ll{x(\gB\log_3x+c_4)^{k-1}\over (k-1)!(\log x)^{A+1}},$$
provided $\K$, and therefore $\kappa$, is suitably chosen. This is the required estimate.
\qed
\medskip\goodbreak
\paradeuxb{Completion of the argument}
It is now a simple matter to derive \eqref{fSD}. Recalling definition \eqref{defTrho}, we have
$$M(x;f)=\sum_{n\leqslant x}g(n)T_\varrho\Big({x\over n}\Big).$$
\par
By \eqref{omsomgj} and partial summation, we have
 $$\gtg_j(y)=\sum_{n\leqslant y}{g(n)(\log n)^j\over n}\ll j(\log y)^{j-A}(\log_2y)^{\gB+ \delta_{A,J+1}}\qquad (y\geqslant 3,\,j>J).\eqdef{majgamj*}$$
  Recall the definition of the coefficients $\alpha_h(\varrho)$ appearing in \eqref{flj} and put $$\nu_h:=\alpha_h(\varrho)/h!\Gamma(\varrho-h)\quad(h\geqslant 0).$$ By \citeplus{Te15}{th. II.5.2}, for a suitable constant $b$, we have $  \nu_h\ll(bh+1)^h$
and $$T_\varrho(x)=x(\log x)^{\varrho-1}\bigg\{\sum_{0\leqslant h\leqslant H}{\nu_h\over (\log x)^h}+O\bigg({(bH+1)^{H+1}\over (\log x)^{H+1}}\bigg)\bigg\}\quad(x\geqslant 2),\eqdef{SDT}$$
uniformly for $H\ll(\log x)^{1/3}$, say. 
\par \goodbreak
Now Dirichlet's hyperbola formula provides 
$$M(x;f)=U+V+O\Big(x(\log x)^{\varrho-1-A}\Big),$$
with 
$$\eqalign{U&:=\sum_{n\leqslant \sqrt{x}}g(n)T_\varrho\Big({x\over n}\Big),\quad
V:=\sum_{d\leqslant \sqrt{x}}\tau_\varrho(d)\sum_{\sqrt{x}<n\leqslant x/d}g(n).\cr}$$
\par 
\goodbreak
The sum $V$  may be treated as an error term: for $d\leqslant \sqrt{x}$, we have, by \eqref{majGj},
$$\sum_{\sqrt{x}<n\leqslant x/d}g(n)\ll {x(\log_2x)^{\gB} \over d(\log x)^{A+1}}\cdot$$ 
By distributing the variable $d$ into intervals $]\e^{-\delta(j+1)}\sqrt{x}, \e^{-\delta j}\sqrt{x}]$ and approximating the corresponding $n$-range by $[\sqrt{x},\sqrt{x}\e^{\delta j}]$ with $\delta:=1/(\log x)^{\K}$ for suitably large $\K$, we obtain
$$\eqalign{V&\ll   \sum_{\di{j\geq 0}{\e^{j\delta}\leqslant \sqrt{x}}} \bigg| \big\{ G_0\big(\e^{\delta j}\sqrt{x}\big)-G_0\big(\sqrt{x}\big) \big\}
\sum_{\e^{-\delta(j+1)}\sqrt{x}<d\leqslant  \e^{-\delta j}\sqrt{x}}\tau_\rho(d)\bigg|+{x \over (\log x)^{A+1}}
\cr
&\ll {x\delta (\log_2x)^{\gB}\over (\log x)^{A+1}} \sum_{\di{j\geq 0}{\e^{j\delta}\leqslant \sqrt{x}}}
\{(\log x)-2j\delta+1\}^{\re \rho-1}
\cr&\ll  {x(\log_2x)^{\gB}\over (\log x)^{A+1}}\Big\{1+ (\log_2x)^{\delta_{0,\re \rho}}(\log x)^{ \re\rho}\Big\} ,
}$$
with Kronecker's notation.
\par 
Finally, we apply \eqref{SDT} with $H=\plaf{A+2r}$ to get
$$U=x\sum_{n\leqslant \sqrt{x}}{  g(n)\over n}\bigg\{\sum_{0\leqslant h\leqslant H}\nu_j  \Big(\log {x\over n}\Big)^{\varrho-h-1}+O\Big({ (\log x)^{\re \varrho-H-2}}\Big)\Big\}.$$
\goodbreak
From the expansion
 $$\bigg({\log x/n\over \log x}\bigg)^{\varrho-h-1}
=
\sum_{0\leqslant k\leqslant K} c_{ kh}\Big({\log n\over \log x}\Big)^k+O\Big(K^{|\rho-1-j|}\Big({\log n\over \log x}\Big)^K\Big)\quad(1\leqslant n\leqslant \sqrt{x}) $$ 
with $c_{kh}\ll k^{|\rho-1-h|},$ we get, for each  $h\in[0,H]$,
$$\eqalign{\sum_{n\leqslant \sqrt{x}}{  g(n)\over n} &\bigg({\log {x/n}\over \log x}\bigg)^{\varrho-h-1}=\sum_{0\leqslant k\leqslant K} c_{kh} {\gtg_k\big(\sqrt{x}\big)\over (\log x)^k} +O\bigg( {K^{|\rho-1-h|}\over 2^{K}}( \log_2  x)^{\gB} \bigg),\cr}
\eqdef{estsum}$$
where we used the trivial estimate $|g|\leqslant \gB^{\omega}$  to bound the error term. 
\par 
Select $K:=\fl{2(A+\gB+1)\log_2x}$. Applying \eqref{omsomgj} for $j\leqslant J$, and \eqref{majgamj*} when $J<j\leqslant H$, we get
$$U=x(\log x)^{\varrho-1}\bigg\{\sum_{0\leqslant j\leqslant J}{\lambda_j(f)\over (\log x)^j}+O\bigg({(\log_2x)^{\gB+\delta_{A,J+1}}\over (\log x)^{A}}\bigg)\bigg\},$$
 the exponent $\delta_{A,J+1}$ arising from that in \eqref{majgamj*}.  This completes the proof.
 \bigskip\medskip
  \paraun{Proof of \ref{SDf}---friable summation}
  \paradeuxb{Setting}
An integer $n$ such that $P^{+}(n)\leqslant y$ is said to be $y$-friable. Friable summability of series, was defined in \citer{Du57}, \citer{FT91}, and has been employed systematically in \citer{Br98}, \citer{BT04}, \citer{BT15}. A series $\sum_{n\geqslant 1}a_n$ is said to be {\it friably summable} to $a$ (or is said to have friable sum $a$) if the subseries  
$$\sum_{\di{n\geqslant 1}{P^+(n)\leqslant y}}a_n$$
converges for each $y\geqslant 2$ and tends to $a$ as $y$ tends to infinity. We then write
$$\sum_{n\geqslant 1}a_n=a\quad (P).$$
\par \goodbreak
Letting $\zeta(s)$ denote Riemann's zeta function, it is well known that, for any given real number $\tau\neq0$, the series $\sum_{n\geqslant 1}1/n^{1+i\tau}$ has friable sum $\zeta(1+i\tau)$  (by the convergence of the Eulerian product on the pointed line $1+i\tau$, $\tau\neq0$), while being divergent in the ordinary meaning.
\par 
We shall show that, representing $f$ as a Dirichlet convolution $f=\tau_\varrho*g$, then, for $0\leqslant j\leqslant J$, we have
$$\gamma_j(g):=\sum_{n\geqslant 1}{g(n)(\log n)^j\over n}\qquad (P).\eqdef{defgamj}$$
With this notation, the coefficients $\lambda_j(f)$ appearing in \eqref{fSDf} are  given by
$$\lambda_j(f)={1\over \Gamma(\varrho-j)}\sum_{\ell+h=j}{\alpha_\ell(\varrho)\gamma_h(g)\over \ell!h!}\qquad (0\leqslant j\leqslant J),\eqdef{flj}$$ while \eqref{flambdaj} remains valid.
\par 
\medskip
\paradeuxb{Reduction}
 Let $g$ be exponentially multiplicative, \ie such that $g(\pnu)=g(p)^\nu/\nu!$ for all primes $p$ and all integers $\nu\geqslant 0$,\note{ This concept has been extensively used in the literature, in particular by Wirsing \citer{Wi67}.}  and define $g(p):=f(p)-\varrho$. Put  $\gf:=\tau_\varrho*g$. Then $f=\gf*h$ where $h$ is multiplicative, supported on the set of squareful integers, and, for each $p$, the values $h(\pnu)$ are given by the power series expansion
$$\sum_{\nu\geqslant 0}h(\pnu)\xi^\nu=(1-\xi)^\varrho\e^{-\xi g(p)}\sum_{\nu\geqslant 0}f(\pnu)\xi^\nu\qquad (|\xi|<1/p^\sigma).$$ 
Thus, $\sum_{n\geqslant 1}|h(n)|/n^\tau<\infty$ for all $\tau>\sigma$ and so, by a standard convolution argument left to the reader, we may restrict to proving \eqref{fSDf} for $\gf$.
\par 
We note right away that  $g\in\FF(2r,\sigma)$ and that \eqref{hypg} holds. 
 Introducing the series
$$\sum_{P^+(n)\leqslant y}{g(n)\over n^s}=\exp\bigg\{\sum_{p\leqslant y}{g(p)\over p^s}\bigg\}\qquad (\re s \geqslant  1),$$
we readily see by partial summation that the right-derivatives  of any order $j\leqslant J$ at  $s=1$ of the left-hand side   converge to a limit as $y\to\infty$, in other words that \eqref{defgamj} holds.  Moreover partial summation yields that,  for  integer $j$ and $y\geqslant 2$,
$$\gamma_j(y;g):=\sum_{P^+(n)\leqslant y}{g(n)(\log n)^j\over n}=\normalbaselineskip=15pt\cases{\gamma_j(g)+O\big((\log y)^{j-A}\big) & \quad $(0\leqslant j\leqslant J),$\cr
O\Big(j!(\log y)^{j-A}(\log_2y)^{\delta_{j,A}}\Big)&\quad $(j>J)$.\cr}\eqdef{somgy}$$
 
 To prove this, observe that, for $|w|\leqslant 1/\log y$, we have
$$\eqalign{\sum_{p\leqslant y}{g(p)\over p^{1-w}}&=\sum_{k\geqslant 0}\sum_{p\leqslant y}{g(p)w^k(\log p)^k\over k! p}=\sum_{k\leqslant J}\bigg\{{\mu_kw^k\over k!}+O\bigg({1\over k!(\log y)^A}\bigg)\bigg\}+O\bigg({(\log_2y)^{\delta_{A,J+1}}\over (\log y)^A }\bigg)\cr
&=\sum_{k\leqslant J}{\mu_kw^k\over k!}+O\bigg({(\log_2y)^{\delta_{A,J+1}}\over (\log y)^A }\bigg),\cr}$$
with $\mu_k:=\sum_{p}g(p)(\log p)^k/p$ $(0\leqslant k\leqslant J)$.
Estimate \eqref{somgy} then follows from Cauchy's formula.  

\medskip
\paradeuxb{An auxiliary estimate}
 For $f=\tau_\varrho*g$ and $g$ as above, let us write $f=\tau_\varrho*g_y*h_y$, where $g_y$ and $h_y$ are the multiplicative functions defined by
 $$g_y(\pnu):=\1_{\{p\leqslant y\}}g(\pnu),\quad h_y(\pnu)=\1_{\{p> y\}}g(\pnu)\qquad (p\geqslant 2,\,\nu\geqslant 1).$$
\par 
Our first goal is to estimate
$$H_y(x):=\sum_{y<n\leqslant x}h_y(n)\qquad (2\leqslant y\leqslant x).$$
 
\Propl{Hy}{Uniformly for  $x\geqslant y\geqslant 2,$ and with $u:=(\log x)/\log y$, we have
$$H_y(x)\ll {xu^{2r} \over (\log x)^{A+1}}\cdot\eqdef{majHy}$$}
\nid 
We have
$$\eqalign{H_y(x)&=\sum_{\di{y<p\leqslant x}{\pnu\leqslant x}}g(\pnu)+\sum_{y<m\leqslant x/y}h_y(m)\sum_{\di{\pnu\leqslant x/m}{p>P^+(m)}}g(\pnu)\ll{x\over (\log x)^{A+1}}+S,\cr}$$
with 
$$\eqalign{S:=
\sum_{\di{y<m\leqslant x/y}{mP^+(m)\leqslant x}}{|h_y(m)|x\over m\{\log (x/m)\}^{A+1}} .\cr}$$
\goodbreak
Writing $m=nq^\nu$ with   $P^+(n) < q$, where, here and in the sequel of this proof, $q$ denotes a prime number, we have 
$$\eqalign{S =\sum_{\di{q>y}{\nu\geqslant 1}}{|g(q)|^{\nu}\over {\nu}!q^{\nu}}\sum_{\di{n\leqslant x/q^{{\nu}+1}}{P^+(n)< q}}{|h_y(n)|x\over n\{\log (x/n q^{\nu})\}^{A+1}},\cr}$$
. \par 
Let $S^-$ denote the contribution to $S$ of $n=1$, and $S^+$ that of $n>y$. We plainly have
$$\eqalign{S^-&\ll \sum_{\nu\geqslant 1}\sum_{y<q\leqslant x^{1/({\nu}+1)}}{|g(q)|^{\nu} x\over {\nu}!q^{\nu}(\log x/q^{\nu})^{A+1}}
\cr&\ll\sum_{{\nu}\geqslant 1}\sum_{y<q\leqslant x^{1/({\nu}+1)}}{x|g(q)|^{\nu} ({\nu}+1)^{A+1}\over {\nu}!q^{\nu}(\log x)^{A+1}}
\ll{x\log (2u)\over (\log x)^{A+1}}\cdot\cr}
$$
\par 
In order to majorize~$S^+$, we appeal to the Rankin-type bound
$$\sum_{\di{y<n\leqslant t}{P^+(n)\leqslant q}}|h_y(n)|\ll t^{1-c/\log q}{(\log t)^{2r-1}\over (\log y)^{2r}}\quad(t, q>y),$$
where $c>0$ is an abolute constant. It follows that
$$\eqalign{S^+&\ll \sum_{{\nu}\geqslant 1}\sum_{y<q\leqslant x^{1/({\nu}+1)}}{x|g(q)|^{\nu}\over {\nu}!q^{\nu}}\int_y^{x/q^{{\nu}+1}}{1\over t(\log x/tq^{\nu})^{A+1}}\d O\Big({t^{1-c/\log q}(\log t)^{2r-1}\over (\log y)^{2r}}\Big)\cr
&\ll \sum_{{\nu}\geqslant 1}\sum_{y<q\leqslant x^{1/({\nu}+1)}}{x|g(q)|^{\nu}\over {\nu}!q^{\nu}(\log y)^{2r}}\int_y^{x/q^{{\nu}+1}}{(\log t)^{2r-1}\d t\over t^{1+c/\log q}(\log x/tq^{\nu})^{A+1}}\cr
&\ll
\sum_{{\nu}\geqslant 1}\sum_{y<q\leqslant x^{1/({\nu}+1)}}{x|g(q)|^{\nu}(\log q)^{2r-A-1}\over {\nu}!q^{\nu}(\log y)^{2r}}\int_0^{(\log x)/(\log q)- {\nu} -1}\!\!\!{\e^{-cw}w^{2r-1}\d w\over \{(\log x)/(\log q)- { {\nu}}-w \}^{A+1}}\cr
&\ll
\sum_{{\nu}\geqslant 1}\sum_{y<q\leqslant x^{1/({\nu}+1)}}{x|g(q)|^{\nu}(\log q^{ {\nu}})^{2r}\over {\nu}!q^{\nu}(\log y)^{2r}(\log x)^{A+1}}\ll{xu^{2r}\over (\log x)^{A+1}}\cdot\cr}$$
 
\vskip-5mm\qed 

\par 
\bigskip
\paradeuxb{Completion of the argument}
 Put $f_y:=\tau_\varrho*g_y$ and let $y:=x^{c /\log_2x }$,  where $c$ is a sufficiently small constant. We have
$$M(x;f_y)=V_y(x)+E_y(x)$$
with 
$$\eqalign{V_y(x)&:=\sum_{n\leqslant \sqrt{x}}g_y(n)T_\varrho\Big({x\over n}\Big),\quad
E_y(x):=\sum_{d\leqslant \sqrt{x}}\tau_\varrho(d)\sum_{\sqrt{x}<n\leqslant x/d}g_y(n).\cr}$$
\par 
We immediately note that a Rankin-type estimate such as \citeplus{Te87}{lemma 2} provides
$$E_y(x)\ll x^{1-1/\log y}\sum_{d\leqslant x}{\tau_r(d)\over d}\ll{x(\log x)^{\varrho-1-N}},\eqdef{evalV''}$$
for any given integer $N$.\par \goodbreak
From \eqref{SDT}, we deduce, for any fixed integer $H$, that
$$V_y(x)=x\sum_{n\leqslant \sqrt{x}}{g_y(n)(\log x/n)^{\varrho-1}\over n}\Bigg\{\sum_{0\leqslant h\leqslant H}{\nu_h\over (\log x/n)^h}+O_H\bigg({1\over (\log x)^{H+1}}\bigg)\Bigg\}.$$
Using the  bound  provided by \eqref{FFB} and selecting $H$ sufficiently large in terms of $r$, we get that the overall contribution of the above remainder term  is
$$\ll x(\log x)^{\varrho+2r-H-2} \ll x(\log x)^{\varrho-1-N}, $$
for any given integer $ N$.
\par Appealing to \eqref{estsum} for $g_y$, we may write
$$\eqalign{\sum_{n\leqslant \sqrt{x}}{  g_y(n)\over n}& \bigg({\log {x/n}\over \log x}\bigg)^{\varrho-h-1}=\sum_{0\leqslant k\leqslant K} c_{kh} \sum_{n\leqslant \sqrt{x}}{g_y(n)\over n} \Big({\log n\over \log x}\Big)^k +O\Big( {1\over  (\log x)^{N+r+1}} \Big),\cr}
 \eqdef{devgy}$$ for sufficiently large $K$.
  Let $v:=1/\log y$. For all $0\leqslant k\leqslant K$, we have, by Rankin's device again, provided $c$ is suitably chosen,
$$\eqalign{\sum_{n>\sqrt{x}}{|g_y(n)|(\log n)^k\over n}&\leqslant {k!\over v^kx^{v/4}}\exp\bigg\{\sum_{p\leqslant y}{|g_y(p)|\over p^{1-v/2}}\bigg\}\ll {k!\over x^{v/4}}(\log y)^{k+2r}\ll{1\over (\log x)^{N+r+1}}\cdot\cr}$$ 
 We may therefore extend, in \eqref{devgy},  the inner $n$-sum  to all $y$-friable integers without perturbing the error. 
Taking \eqref{evalV''} into account, we thus obtain
$$M(x;f_y)=x(\log x)^{\varrho-1} \bigg\{\sum_{0\leqslant h\leqslant H}{\nu_h\over (\log x)^h}\sum_{ 0\leqslant k\leqslant K}c_{kh}{\gamma_k(y;g)\over (\log x)^k}+O\bigg({1\over (\log x)^N}\bigg)\bigg\},$$
 which, by rearranging the double sum, yields
$$M(x;f_y)=x(\log x)^{\varrho-1}\bigg\{\sum_{0\leqslant j\leqslant  K+H}{\lambda_j(y;f)\over (\log x)^j}+O\bigg({1\over (\log x)^{N}}\bigg)\bigg\}\eqdef{evalVy}$$
where, in view of \eqref{somgy}, $$\lambda_j(y;f)=\lambda_j(f)+O\Big(j^j(\log y)^{j-A}(\log_2y)^{\delta_{A,j}}\Big)\qquad (0\leqslant j\leqslant K+H),\eqdef{majlambdaj}$$
with the convention that, say, $\lambda_j(f)=0$ for $j>J$.
\par \goodbreak
By Dirichlet's hyperbola formula, we have
$$M(x;f-f_y)=M(x;f_y*h_y-f_y)=S+U-W,\eqdef{sum-sumfy}$$
with
$$S:=\sum_{n\leqslant \sqrt{x}}f_y(n)H_y\Big({x\over n}\Big),\quad U:=\sum_{y<m\leqslant \sqrt{x}}h_y(m)M\Big({x\over m};f_y\Big),\qquad W:=M\big(\sqrt{x};f_y\big)H_y\big(\sqrt{x}\big).$$
 \par 
 Now
$$S\ll\sum_{n\leqslant \sqrt{x}}{|f_y(n)|x(\log_2x)^{2r}\over n(\log x)^{A+1}}\ll {x(\log_2x)^{2r}\over (\log x)^{A+1-r}},$$
and similarly
$$W\ll x(\log x)^{\varrho-A-2}(\log_2x)^{2r}
\ll {x(\log_2x)^{2r}\over (\log x)^{A+1-r}}.$$ 
\par \goodbreak
 To estimate the sum $U$, we insert \eqref{evalVy} and evaluate the contribution of the main terms by partial summation from \eqref{majHy}. Taking
\eqref{majlambdaj} into acount, we obtain 
$$U\ll x u^{2r}(\log x)^{\varrho-A-1}.$$
 \par 
Recalling \eqref{evalVy}, \eqref{majlambdaj} and \eqref{sum-sumfy}, we have thus almost reached  \eqref{fSDf}, but with a remainder term
 $$\R(x;f)\ll {x(\log_2x)^{2r}\over(\log x)^{A+1-r} }\cdot$$  
\par 
We now show that an inductive argument enables us to replace $(\log_2x)^{2r}$ by a power of~$\log_3x$. 
Indeed, under assumption \eqref{FFB}, consider the exponentially multiplicative function~$\varphi$ defined~by
$$\varphi(\pnu)=f(p)^\nu/2^\nu\nu!\quad(p\geqslant 2,\,\nu\geqslant 1).$$
 Then $\varphi\in\FF(r/2,\sigma)$. Moreover, we may write $f=\varphi*\varphi*\psi$ with
$$\psi(\pnu)=\sum_{0\leqslant j\leqslant\nu}{(-1)^jf(p)^jf(p^{\nu-j})\over j!}\quad(\nu\geqslant 1).$$
Thus $\psi(p)=0$ and, writing $\varepsilon_p:=|f(p)|^2/p^{2\sigma}+\sum_{\nu\geqslant 2}|f(\pnu)|/p^{\nu\sigma}$, we have
$$\eqalign{\sum_{\nu\geqslant 2}{|\psi(\pnu)|\over p^{\sigma\nu}}&\leqslant \sum_{j+k\geqslant 2}{|f(p)|^j|f(p^k)|\over j!p^{\sigma(k+j)}}\cr
&\leqslant \varepsilon_p+{|f(p)|\over p^{\sigma}}\varepsilon_p+\varepsilon_p\e^{|f(p)|/p^{\sigma}}\bigg\{1+{|f(p)|\over p^\sigma}+\varepsilon_p\bigg\}\ll\varepsilon_p.\cr}$$
Applying the estimate already proved for $\varphi$ and writing the hyperbola formula for $f$ furnishes
$$\R(x;f)\ll{x(\log_2x)^{r}\over(\log x)^{A+1-r} }\cdot$$
 After $k+1$ iterations, the exponent of $\log_2 x$ becomes $r/2^k$, while the implicit constant in the upper bound for $\R(x;f)$ gets  multiplied by $C^k$ for a suitable constant $C$. Selecting $k\asymp \log_4x$ yields \eqref{fSDf}. 
\bigskip\medskip
\paraun{Limitations}
\paradeuxb{Optimatility of \ref{SD+}}
\drefdeux{contrex1}  We show here that it is not possible, in general, to replace $\max(0,\re\varrho)$ by $\re\varrho$ in the error term of \eqref{fSD}.
  To this end, we initially consider the completely multiplicative function $g_0$ defined by $$g_0(\pnu):=1/(\log p)^{A\nu}\quad(p\geqslant 2,\,\nu\geqslant 1)$$ with $0<A<1$---note however that, at the cost of slightly more complicated computations, our approach would work for any positive $A$. 
As a preliminary step we show that, for a suitable  constant $C$,  we have
  $$M(x;g_0)={Cx\over (\log x)^{A+1}}\Big\{1+O\Big({1\over (\log x)^{A}}\Big)\Big\}\quad(x\geqslant 2).\eqdef{faG0}$$
 Indeed, on the one hand, since $\sum g_0(n)/n<\infty$, we have
 $$M(x;g_0)\ll x/\log x\quad(x\geqslant 2),\eqdef{maj1G}$$
 and, on the other hand,
$$\eqalign{M(x;g_0)\log x-&\int_1^x{M(t;g_0)\over t}\d t=\sum_{m\leqslant x}g_0(m)\sum_{d\leqslant x/m}g_0(d)\Lambda(d)\cr &=x\sum_{m\leqslant x}{g_0(m)\over m(\log 2x/m)^A}+O\Big(x\sum_{m\leqslant x}{g_0(m)\over m(\log 2x/m)^{A+1}}\Big).\cr}\eqdef{est1G}$$
From this and \eqref{maj1G}, we deduce that $M(x;g_0)\ll{x/(\log x)^{A+1}}$ $(x\geqslant 2)$. Carrying this estimate back into~\eqref{est1G} yields
$$\eqalign{M(x;g_0)&={x\over \log x}\sum_{m\leqslant x}{g_0(m)\over m(\log 2x/m)^A}+O\Big({x\over (\log x)^{A+2}}\Big).\cr
&={x\over \log x}\sum_{m\leqslant \sqrt{x}}{g_0(m)\over m   (\log x/m )^A } 
 +O\Big( {x\over (\log x)^{ 2A+1}}\Big)\cr
&={x\over  (\log x )^{A+1}}\sum_{m\leqslant \sqrt{x}}{g_0(m)\over m }\Big\{1+O\Big( {\log m\over  \log x }\Big)\Big\} +O\Big( {x\over (\log x)^{ 2A+1}}\Big)
\cr
&={x\over  (\log x )^{A+1}}\sum_{m\leqslant \sqrt{x}}{g_0(m)\over m } 
 +O\Big( {x\over (\log x)^{ 2A+1 }}\Big),\cr} $$
which implies \eqref{faG0} with
$C:=\sum_{m\geqslant 1}{g_0(m)/ m } .$
 \goodbreak
 \par 
Next, define $g(n):=g_0(n)n^i$ $(n\geqslant 1)$. By partial summation, we have
$$M(x;g)={Cx^{1+i}\over (1+i)(\log x)^{A+1}}\Big\{1+O\Big({1\over (\log x)^{A}}\Big)\Big\}\quad(x\geqslant 2).\eqdef{faG}$$
Now,  select $\varrho=-r$, with $A<r<1$ and define $f=g*\tau_\varrho$. This function satisfies the hypotheses of \ref{SD+}. By the hyperbola formula, we have, for $1\leqslant y\leqslant x$,
$$M(x;f)=U+V-W\eqdef{hypf}$$
with 
$$U:=\sum_{n\leqslant x/y}g(n)T_\varrho\Big({x\over n}\Big),\quad V:=\sum_{n\leqslant y}{\tau_\varrho(n)}M\Big({x\over n};g\Big),\quad W:=T_\varrho(y)M\Big({x\over y};g\Big).$$
We choose $y:=\e^{a(\log_2x)^2}$, where $a$ is a sufficiently large constant. We plainly have 
$$W\ll {x\over (\log y)^{r+1}(\log x)^{A+1}}.$$
The sum $U$ may be evaluated by the Selberg-Delange method in the form given in \citeplus{Te15}{\S\thinspace II.5.4, Notes}. For suitable positive constants  $b$, $c$, we have
$$T_\varrho(x)=\int_0^b\alpha(t)x^{1-t}t^r\d t+O\Big(x\e^{-c\sqrt{\log x}}\Big)\qquad (x\geqslant 2),$$
where $\alpha$ is continuous on $[0,b]$. Since we may deduce from \eqref{faG} that
$$Z(t):=\sum_{n\leqslant x/y}{g(n)\over n^{1-t}}\ll{(x/y)^t\over (\log x)^{A+1}}\qquad (0\leqslant t\leqslant b),$$
we get, with an appropriate choice of $a$, 
$$U=\int_0^b\alpha(t)x^{1-t}Z(t)t^r\d t+O\Big({x\over (\log x)^{A+2}}\Big)\ll{x\over (\log x)^{A+1}(\log y)^{r+1}}\cdot$$
Finally
$$\eqalign{V&={Cx^{1+i}\over (1+i)(\log x)^{A+1}}\sum_{n\leqslant y}{\tau_\varrho(n)\over n^{1+i}}\Big\{1+O\Big({\log n\over \log x}\Big)\Big\}+O\bigg({x\over (\log x)^{2A+1}}\sum_{n\leqslant y}{\tau_r(n)\over n}\bigg)\cr
&= {\{1+o(1)\}Cx^{1+i}\over (1+i)\zeta(1+i)^r(\log x)^{A+1}}\cdot\cr}$$
Carrying back into \eqref{hypf}, we obtain
$$M(x;f)\gg{x\over (\log x)^{A+1}},$$
which implies the stated property.
\medskip\bigskip
\paradeuxb{Optimality of \ref{SDf}}
\drefdeux{contrex2}
We show here that the exponent $r$ appearing in the remainder term of \eqref{fSDf} cannot in general be replaced by $\re\varrho$, even when $\re\varrho\geqslant 0$.  To this end we exhibit a function $f$ belonging to $\FF(\sigma,r)$ for $1/2<\sigma<1$ and for which: (i) condition \eqref{hyp} fails, (ii) $\varrho=0$, $0<A<r$, and (iii)  given any $c<1$, one cannot replace $r$ by $cr$  in \eqref{fSDf}. If $f$ happens to be real, the conditions become $0<A<2r/\pi$ and $c<2/\pi$.
Our construction relies on exploiting possible resonance between $f(p)$ and $p^{it}$ for certain values of the real parameter $t$. 
\par
 Let $A>0$, $r>0$, $C>0$ (large),  and let $x_k:=\exp\exp C^k$, \hbox{$t_k:=(\log x_{k+1})^A$~$(k\geqslant 1)$}. Put $\varphi(v):=\e^{iv}$  and  define $f$ as the exponentially multiplicative function such that \hbox{$f(p):=r\varphi(t_k\log p)$} whenever $x_k<p\leqslant x_{k+1}$ and, say, $f(p)=0$ when $p\leqslant x_1$. Observe that $\int_0^y\varphi(tw)\d w\ll1/t$ $(t\geqslant 1,\,y>0)$. Due to the rapid increase of $x_k$, partial summation hence yields
$$\sum_{p\leqslant x}f(p)\log p\ll{x\over (\log x)^A}\qquad (x\geqslant 2).$$
Now, for $x:=x_{k+1}$, $t:=t_k$, $\sigma:=1+1/\log x$, we have
$$\abs{\sum_{n\geqslant 1}{f(n)\over n^{\sigma+it}}}\asymp\prod_{h\leqslant k}\e^{S_h}\ll t\int_1^x{|M(y;f)|\over y^2}\d y\eqdef{minmoyf}$$
with $$S_h:=r\sum_{x_h<p\leqslant x_{h+1}}{\re\{f(p)/p^{it}\}\over p}\qquad (1\leqslant h\leqslant k).$$
\goodbreak
Of course $$\eqalign{S_k:&=r{\log \Big({\log x_{k+1}\over \log x_k}\Big)}+O(1)= rC^{k+1}(1-1/C)+O(1).\cr}\eqdef{estSk}$$
Next, for $2\log k<h\leqslant k-1$, we have, by the prime number theorem in a strong form,
$$\eqalign{S_h&=r\int_{x_h}^{x_{h+1}}{\cos\{(t_k-t_h)\log v)\}\over v\log v}\d v+O(1)\cr
&=r\int_{\log x_h}^{\log x_{h+1}}{\cos\{(t_k-t_h)w\}\over w}\d w+O(1)\ll{1\over t_k\log x_h}+1\ll1.\cr}$$
Bounding  $S_h$ trivially for $h\leqslant 2\log k$, we finally get
$$\sum_{h\leqslant k}S_h\geqslant r C^{k+1}(1-1/C)+O\big(k^{2\log C}\big)= r(1-1/C) \log_2x+O\Big((\log_3x)^{2\log C}\Big).$$
\goodbreak
Carrying back into \eqref{minmoyf} yields
$$\int_1^x{|M(y;f)|\over y^2}\d y\gg(\log x)^{r(1-1/C)-A+o(1)}.$$
and so, as $x\to\infty$, provided $A<r(1-1/C)$,
$$M(x;f)=\Omega\bigg({x\over (\log x)^{A+1-r(1-1/C)+o(1)}}\bigg).\eqdef{minMfx}$$
Since $C$ may be chosen arbitrarily large, this furnishes the required result.
\par \goodbreak
 If we require $f$ to be real, we define $\varphi(v):=\sgn(\cos v)$ $(v\in\r)$.  Since $\varphi(v)\cos v$ has mean value $2/\pi,$ we obtain $S_k=(2r/\pi)C^{k+1}(1-1/C)+O(1)$ instead of \eqref{estSk}, appealling, \eg, to \citeplus{Te15}{Lemma III.4.13}. The remainder of the analysis is essentially identical and so we get~\eqref{minMfx} with $2r/\pi$ in place of $r$. It is noticeable that, provided  $A$ is suitably chosen, the implied lower bound  and \eqref{majHT+} may
agree to an arbitrary small power of $\log x$. 
\par \medskip
 \rem The above construction shows that assumptions $f\in\FF(\sigma,r)$ $(\sigma<1,\, r>1)$ and \eqref{hypg} with $\varrho=0$, $0<A<r-1$, are insufficient to imply the convergence of the series $\sum_{n\geqslant 1}f(n)/n$. 
\bigskip
\goodbreak 
\centerline{\twelvebf References}
\medskip
\eightpoint{
\bibtem{Br98} R. de la Bretche, $P$-rŽgularitŽ de sommes d'exponentielles, {\it Mathematika}, {\bf 45} (1998), 145--175.\par
\bibtem{BT04} R. de la Bretche \& G. Tenenbaum, SŽries trigonomŽtriques ˆ coefficients arithmŽtiques,
{\it J. Anal. Math. \bf92} (2004), 1--79.\par
\bibtem{BT15} R. de la Bretche \& G. Tenenbaum, ThŽorme de Jordan friable, in: {\it Analytic number theory, volume in honor of Helmut Maier}, C. Pomerance, M. Rassias, eds., Springer, New York 2015, 57--64. \par 
\bibtem{De59} H. Delange, Sur des formules dues \`a Atle Selberg, {\it Bull. Sc. Math. (2)
\bf 83} (1959), 101--111.
\bibtem{De71} H. Delange, Sur des formules de Atle Selberg, {\it Acta Arith. \bf 19} (1971), 105--146.
\bibtem{Du57}   R.J. Duffin, Representation of Fourier integrals as sums, III. {\it Proc. Amer. Math. Soc. \bf8} (1957), 272--277.\par
\bibtem{FT91} ƒ. Fouvry \& G. Tenenbaum, Entiers sans grand facteur premier en
progressions arithmŽtiques, {\it Proc. London  Math. Soc.} (3) {\bf 63} (1991), 449--494.\par
\bibtem{GK19} A. Granville \& D. Koukoulopoulos, Beyond the LSD method for the partial sums of multiplicative functions, {\it Ramanujan J. \bf49} (2019), 287-319.\par 
\bibtem{Ha95} R.R. Hall, A sharp inequality of Hal‡sz type for the mean value of a multiplicative arithmetic function, {\it  Mathematika \bf42}, \numero1 (1995), 144--157.\par 
\bibtem{HT91}  R.R. Hall \& G. Tenenbaum, Effective mean value estimates for complex multiplicative
functions, {\it Math. Proc. Camb. Phil. Soc.} {\bf 110} (1991), 337--351.\par
\bibtem{HTW08} G. Hanrot, G. Tenenbaum \& J. Wu, Moyennes de certaines fonctions
multiplicatives sur les entiers friables, 2, \PLMS\ (3) \bf96 \rm(2008), 107--135.
\bibtem{Sa53} L.G. Sathe, On a problem of Hardy on the distribution of integers having a
given number of prime factors, I, II, {\it J. Indian Math. Soc. \bf 17} (1953), 63--141.
\bibtem{Sa54} L.G. Sathe, On a problem of Hardy on the distribution of integers having a
given number of prime factors, III, IV, {\it J. Indian Math. Soc. \bf 18} (1954), 27--81.
\par 
\bibtem{Se54} A. Selberg, Note on the paper by L.G. Sathe, {\it J. Indian Math. Soc. \bf 18} (1954),
83--87. 
\bibtem{Sh80} P. Shiu, A Brun--Titchmarsh theorem for multiplicative functions, {\it J. reine
angew. Math. \bf313} (1980), 161--170.\par 
\bibtem{Te87} G. Tenenbaum, Un problme de probabilitŽ conditionnelle en ArithmŽtique, {\it Acta
Arith. \bf 49} (1987), 165--187.\par 
\bibtem{Te00} G. Tenenbaum, A rate estimate in Billingsley's theorem for the size distribution of large prime
factors, {\it Quart. J. Math. (Oxford) \bf51} (2000), 385--403.
\bibtem{Te15} G. Tenenbaum, {\it Introduction to analytic and probabilistic number theory}, 3rd ed., Graduate Studies in Mathematics 163, Amer. Math. Soc. 2015.\par
\bibtem{Te17} G. Tenenbaum, Moyennes effectives de fonctions multiplicatives complexes, {\it Ramanujan J. \bf44}, \numero3 (2017), 641--701;  Corrig. {\it ibid. \bf 51}, \numero1 (2020), 243-244.\par
\bibtem{TW08} G. Tenenbaum \& J. Wu, Moyennes de certaines fonctions
multiplicatives sur les entiers friables, 3, {\it Compositio Math. \bf144} \numero2 (2008), 339--376. 
\bibtem{Wi61} E. Wirsing, Das asymptotische Verhalten von Summen Ÿber multiplikative Funktionen, {\it Math. Ann. \bf143} (1961), 75Ð102.\par 
\bibtem{Wi67} E. Wirsing, Das asymptotische Verhalten von Summen \"uber multiplikative Funktionen,
II, {\it Acta Math. Acad. Sci. Hung. \bf 18} (1967), 411--467.
\par }
\smallskip
{\leftskip2mm\rightskip-2cm\sevenrm
\gutter=4cm \doublecolumns
 \obeylines \baselineskip=7pt
RŽgis de la Bretche
UniversitŽ de Paris, Sorbonne UniversitŽ, CNRS, 
Institut de Math. de Jussieu-Paris Rive Gauche,
 F-75013 Paris,  
France
\smallskip
{\seventt regis.de-la-breteche@imj-prg.fr}
 GŽrald Tenenbaum\par
Institut \'Elie Cartan\par 
Universit\'e de Lorraine\par
 BP 70239\par
54506 Vand\oe uvre-ls-Nancy Cedex\par
 France
\smallskip
{\seventt gerald.tenenbaum@univ-lorraine.fr}
\singlecolumn
\par}
\end